\documentclass[10pt]{article}

\usepackage{amsmath,amsmath,amssymb}
\usepackage{amsthm}

\usepackage{mathrsfs}
\usepackage{multirow}
\usepackage{graphicx}
\usepackage{authblk}
\usepackage{indentfirst}
\usepackage{multicol}  
\usepackage{tabu}
\usepackage{tabularx}
\usepackage{url}
\usepackage{fancyhdr}
\usepackage[numbers]{natbib}
\usepackage{lineno}

\usepackage{fancyhdr}
\usepackage{colortbl}
\usepackage{caption}

\usepackage{hyperref}
\hypersetup{colorlinks=true,linkcolor=blue,anchorcolor=blue,citecolor=blue}

\usepackage[table,xcdraw]{xcolor}
\usepackage{longtable}
\usepackage{hhline}
\usepackage{alltt}
\usepackage{tgcursor}
\usepackage{xr}

\newtheorem{definition}{Definition}
\newtheorem{proposition}{Proposition}
\newtheorem{theorem}{Theorem}
\newtheorem{example}{Example}
\newtheorem{remark}{Remark}
\makeatletter
\newcommand*{\addFileDependency}[1]{
	\typeout{(#1)}
	\@addtofilelist{#1}
	\IfFileExists{#1}{}{\typeout{No file #1.}}
}
\makeatother
\newcommand*{\myexternaldocument}[1]{
	\externaldocument{#1}
	\addFileDependency{#1.tex}
	\addFileDependency{#1.aux}
}
\myexternaldocument{supplement}

\begin{document}
	\setlength{\parindent}{2em}
	\title{Knot data analysis using multiscale Gauss link integral }

	\author[1]{Li Shen } 
	\author[1]{Hongsong Feng }
  \author[2]{Fengling Li }
	\author[2]{Fengchun Lei }
  \author[3]{Jie Wu } 
	\author[1,4,5]{Guo-Wei Wei \thanks{Corresponding author: weig@msu.edu}}
	\affil[1]{Department of Mathematics, Michigan State University, East Lansing, MI 48824, USA}
  \affil[2]{School of Mathematical Sciences, Dalian University of Technology, Dalian 116024, China }
	\affil[3]{Yanqi Lake Beijing Institute of Mathematical Sciences and Applications, 101408,
China}
	\affil[4]{Department of Biochemistry and Molecular Biology, Michigan State University, MI, 48824, USA}
	\affil[5]{Department of Electrical and Computer Engineering, Michigan State University, MI 48824, USA}
	\renewcommand*{\Affilfont}{\small\it}
	\renewcommand\Authands{ and }
	\date{}
	
	
	\maketitle
	\begin{abstract}
In the past decade, topological data analysis (TDA) has emerged as a powerful approach in data science. The main technique in TDA is persistent homology, which tracks topological invariants over the filtration of point cloud data using algebraic topology. Although knot theory and related subjects are a focus of study in mathematics, their success in practical applications is quite limited due to the lack of localization and quantization. We address these challenges by introducing knot data analysis (KDA), a new paradigm that incorporating curve segmentation and multiscale analysis into the Gauss link integral.  The resulting multiscale Gauss link integral (mGLI) recovers the global topological properties of knots and links at an appropriate scale but offers multiscale feature vectors to capture the local structures and connectivities of each curve segment at various scales. The proposed mGLI significantly outperforms other state-of-the-art methods in benchmark protein flexibility analysis, including earlier persistent homology-based methods. Our approach enables the integration of artificial intelligence (AI) and KDA for general curve-like objects and data.

	\end{abstract}	
	

	\begin{verbatim}
	\end{verbatim}
 {\setcounter{tocdepth}{4} \tableofcontents}
	
	\newpage
	
	\section{Introduction}
	Knot theory is a branch of mathematics that study mathematical knots \cite{crowell2012introduction}. These mathematics knots are often defined as embeddings of a closed circle $S^1$ into the three-dimensional (3D) Euclidean space.  One of the central problems in knot theory is the knot classification problem \cite{adams1994knot}, which asks whether two knots are equivalent or not. Knots that are equivalent can be transformed into each other by a series of continuous deformations, such as stretching or bending the curve, without cutting or gluing any part of it. Such deformations are also known as ambient isotopy \cite{armstrong2013basic}. Generally, mathematicians are accustomed to studying a knot in the sense of ambient isotopy as it allows them to focus on the most essential properties of knots themselves rather than on a particular way in which the knot is embedded in the 3D space. Mathematicians have proposed various knot invariants under ambient isotopy which can classify and characterize different type of knots, such as knot crossing number, knot group \cite{crowell2012introduction}, knot polynomials\cite{adams1994knot}, knot Floer homology \cite{manolescu2014introduction}, Khovanov homology \cite{khovanov2000categorification}, etc.
	
	Knot theory has multiple applications in various fields including physics \cite{ohtsuki2002quantum},  biochemistry \cite{liang1994knots}, and biology \cite{sumners2020role,schlick2021knot,millett2013identifying}. In the above application scenarios, the studied object does not always meet the ideal state expected mathematically. First of all, such target might not be a closed circle. Secondly, the ambient isotopy may lead to major changes in the properties of the object. For instance, the realization of many application object functions depends on local structural information. An ambient isotopy may completely alter the local structure while keeping the global knot information unchanged. Therefore, it is imperative to develop  knot theory-based tools that are robust and effective for data science.
	
	Several attempts have been made to address the aforementioned challenge. Jamroz et al. proposed the protein topology database KnotProt to study  knot and slipknot type of proteins \cite{jamroz2015knotprot}. Dabrowski-Tumanski et al.'s Topoly extends the work to include links and spatial graphs, and also enables the calculation of topological polynomials invariant of those structures \cite{dabrowski2021topoly}.  Recently, Panagiotou and Kauffman have proposed   new invariants for open curves in 3-space \cite{panagiotou2020knot}. Such method  offers a well-defined measurement of entanglement of an open curve. In addition, Baldwin et al. \cite{baldwin2022local} have also made some attempts to localize knot information: intercepting some specific intervals in the linear structure of an open curve gives smaller open curves. Nevertheless, this approach does not really describe local information in 3D space and thus has very limited power in applications related to local functions. In general, current knot and related theories place a major emphasis on preserving global topological properties, and thus, like the classical homology, have very limited success in real-world applications.      
	
We believe that a feasible localization scheme for knot and related theories is to undertake a multiscale analysis. Multisalce analysis has achieved great success in various areas  of mathematics, including wavelet theory, differential equations, and  topological data analysis (TDA). For example, the key technique in TDA is persistent homology, a mathematical framework that combines concepts from algebraic topology,  geometry, and multiscale analysis to analyze complex datasets \cite{edelsbrunner2008persistent,zomorodian2004computing}. It aims to uncover the underlying topological structure and shape of data at different scales, enabling insights that may not be easily discernible with traditional geometric and statistical techniques.  By leveraging the power of topological deep learning (TDL) \cite{cang2017topologynet},  persistent homology offers a unique perspective on understanding the shape and structure of complex datasets, providing valuable insights into the underlying patterns and relationships within the data. TDA has had tremendous success in data science and biological science \cite{chen2022persistent}. 

The linking number is a numerical topological invariant that measures the extent of linkage between two closed curves in three-dimensional space, representing the number of times that each curve winds around the other. The Gauss linking integral \cite{gauss1833integral}, also known as Gauss's integral for the linking number, gives an explicit formulation for the linking number. It is a fundamental tool in the study of knots, links, and other topological structures in three-dimensional space. Gauss linking integral plays an important role in knot theory, algebraic topology,   differential geometry, and quantum field theory. For example, for idealized Dirac-string center vortices,  the Chern-Simons  number can   be given by the Gauss link integral \cite{cornwall2002sphalerons}.  High-order link integrals were proposed \cite{berger1990third}.   However, these approaches are typically global and qualitative.   

The objective of this work is to introduce knot data analysis (KDA) as a new paradigm for data science.  To this end,  we formulate a new framework called  multiscale Gauss linking integral (mGLI) by integrating multiscale analysis with classical knot and knot-related theories.  The proposed mGLI can capture both local and global information of knots, curves, and other curve-like objects. Our framework enables us to analyze curve-like objects by applying a specific radius of an open ball around each segment of curve-like objects and provides a metric that describes the degree of local entanglement. As one gradually increases the radius of the open ball, the newly incorporated region within it will reflect its impact on the original structure by altering the value of the metric. When the radius reaches its maximum limit, the metric will retain the global  information, such as knots and entangled links. To demonstrate the utility of the proposed mGLI, we apply mGLI to protein flexibility analysis, which successfully constructs a robust correlation  between metrics of mGLI and protein   B-factors.  We show that our KDA tool significantly  outperforms other state-of-art methods, including persistent homology, for protein flexibility prediction. 
	
	\section{Methods}
	\subsection{Gauss linking integral of  open and closed curves}

	\begin{definition}[Gauss linking integral]
	Given two disjoint  open or closed curves $l_1 $ and $ l_2$, parametrized as $\gamma_1(s) $ and $ \gamma_2(t)$, respectively. Then the following double integral gives the the Gauss linking integral that  characterizes the degree of interlinking between $l_1$ and $l_2$ \cite{ricca2011gauss}:
	\begin{equation}
		L(l_1,l_2) = \frac{1}{4\pi}	\int_{[0,1]}\int_{[0,1]}\frac{\det(\dot{\gamma_1}(s),\dot{\gamma_2}(t),\gamma_1(s)-\gamma_2(t))}{|\gamma_1(s)-\gamma_2(t)|^3}\ ds\ dt,
	\end{equation}
where $\dot\gamma_1(s)$ and $\dot\gamma_2(t)$ are derivative of $\gamma_1(s)$ and $\gamma_2(t)$, respectively.
	\end{definition}

\begin{proposition}
	The Gauss linking integral is identical to the average of half the algebraic sum of inter-crossings in the projection of the two curves in any possible projection direction for both open and closed curves.
\end{proposition}

\begin{theorem}[Panagiotou et al.\cite{panagiotou2020knot}]
For closed curves, the Gauss linking integral is an integer and a topological invariant. For open curves, the Gauss linking integral is a real number and a continuous function of curve coordinates.
\end{theorem}

\begin{definition}[Segmentation of Gauss linking integral]\label{discrimination}
		Giving finite curve segments  $P_n $ and $ Q_m$ for disjoint open or closed curves $l_1$and  $l_2$, respectively, then the segmentation of Gauss linking integral  induced by the curve segments is defined as the following $n\times m$ segmentation matrix:
		\begin{equation} \label{G}
			G = \begin{pmatrix}
				L(p_1,q_1)&L(p_1,q_2)& \cdots &L(p_1,q_m) \\
				L(p_2,q_1)&L(p_2,q_2)& \cdots &L(p_2,q_m)\\
				\vdots       & \vdots        &\ddots & \vdots \\
				L(p_n,q_1)& L(p_n,q_2)&\cdots&L(p_n,q_m)\\			
			\end{pmatrix},
		\end{equation}
	where $p_i \in P_n$ and $q_j \in Q_m$ are curve segments of $l_1$ and $l_2$, respectively. 
	\end{definition}

\begin{figure}[htb!]
	\centering
	\begin{tabular}{cc}
				\includegraphics[width=0.2\textwidth]{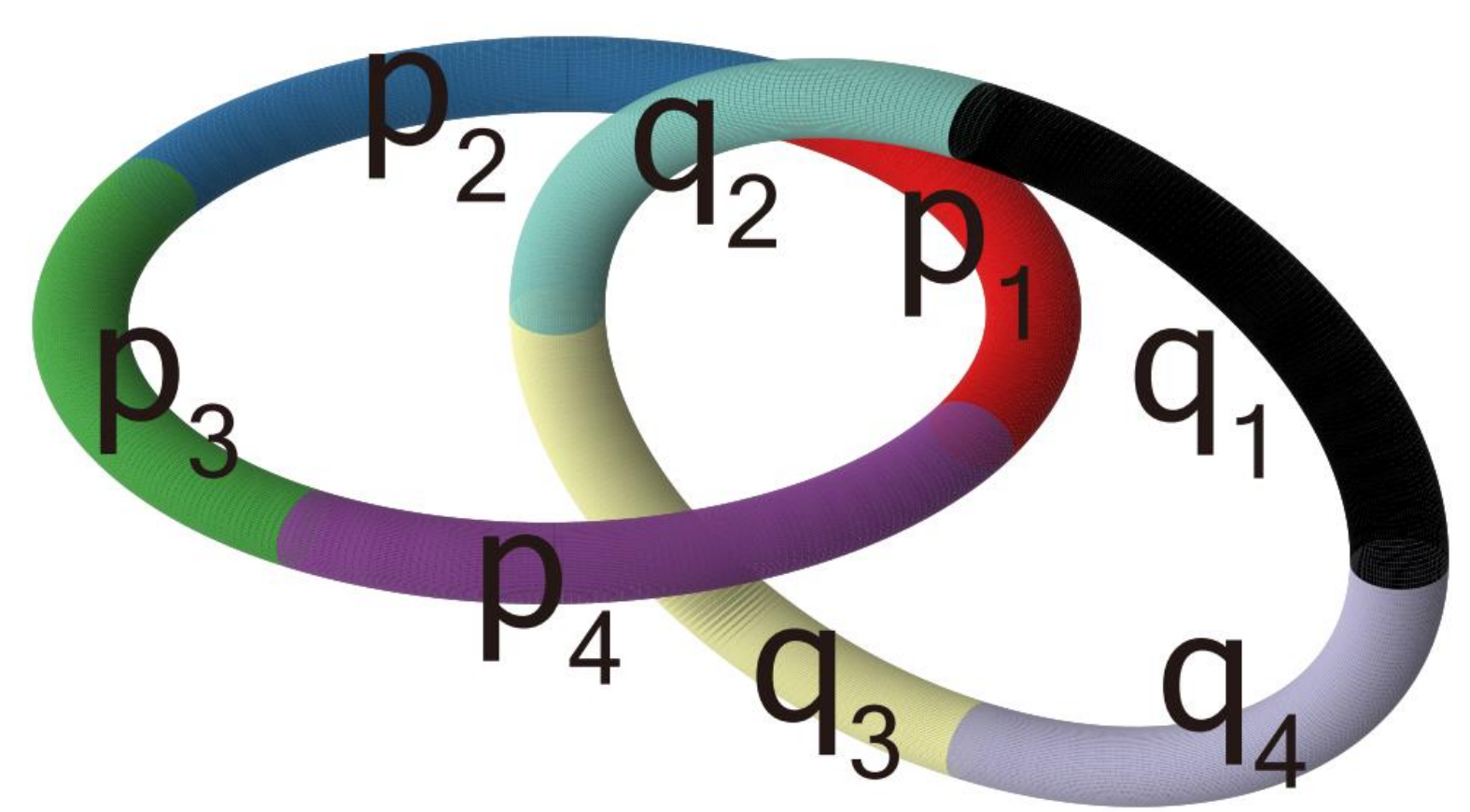}&
				\includegraphics[width=0.2\textwidth]{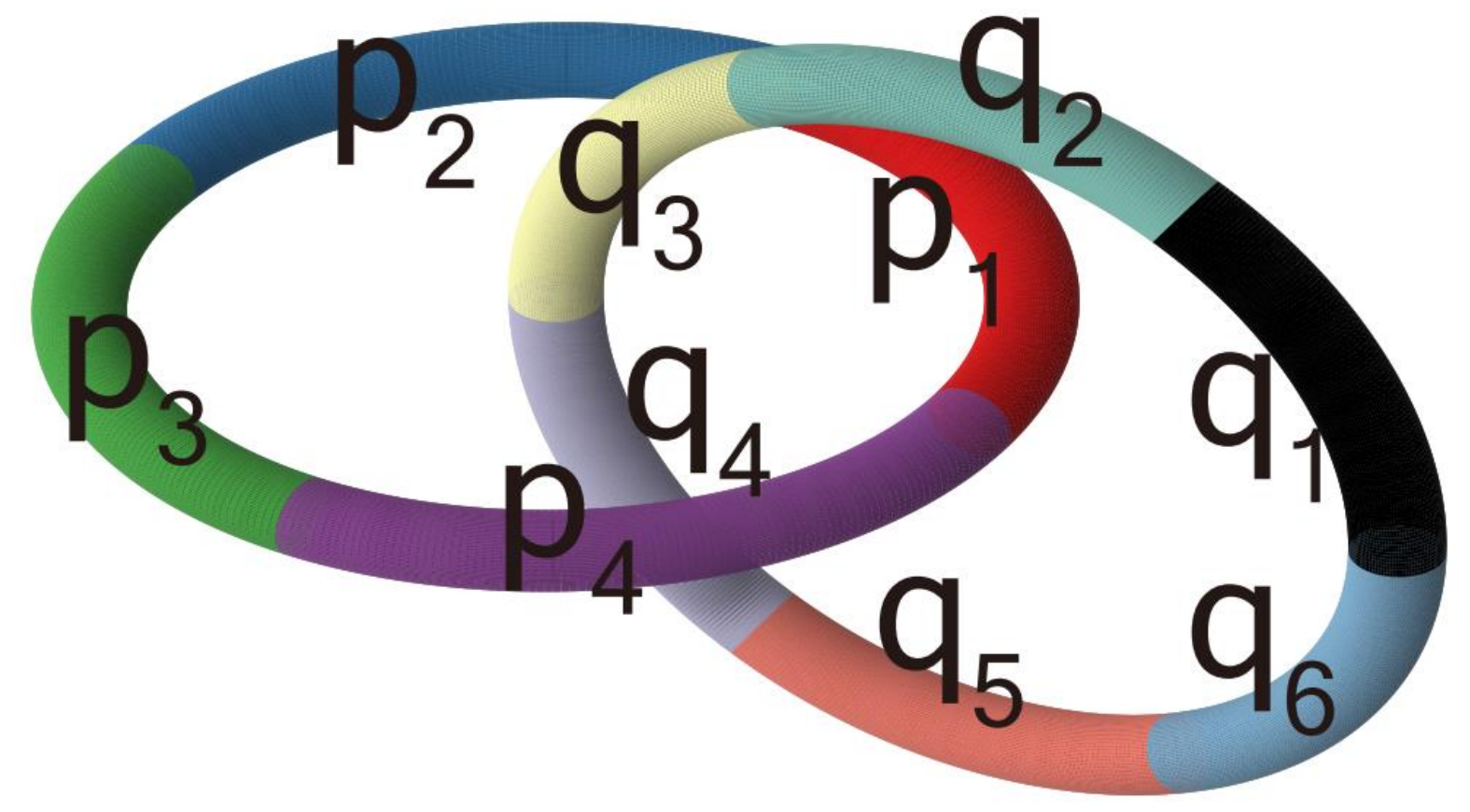}
	\end{tabular}
		\caption{Curve segmentation of the Hopf link. }
\label{example}
\end{figure}	

\begin{theorem}[The grand sum of the segmentation matrix]
	The grand sum of the segmentation matrix of  two curves equals  the Gauss linking integral of the original curves:
	\begin{equation}
		L(l_1,l_2) = \sum_{i}\sum_{j}L(p_i,p_j).
	\end{equation}
\end{theorem}

\begin{example}
\autoref{example} shows two different curve segmentation  of the Hopf  link of two components $l_1$ and $l_2$, with the following parametric equation, respectively:
\begin{equation}
	\gamma_1 = \begin{pmatrix}
		\cos(2\pi t)\\
		\sin(2 \pi t)\\
		0
	\end {pmatrix} ~ {\rm and~ } \gamma_2 = \begin{pmatrix}
	\cos(2\pi s)+1\\
	0\\
	\sin(2\pi s)\
\end{pmatrix}.
\end{equation}

It is well-known that the Gauss linking integral of the Hopf link shown in \autoref{example}  equals to $-1$. Based on \autoref{G}, we can obtained the following segmentation  of Gauss linking integral:
\begin{equation}
		G_1 =  \begin{pmatrix}
 	L(p_1,q_1)&L(p_1,q_2)& L(p_1,q_3) &L(p_1,q_4) \\
 	L(p_2,q_1)&L(p_2,q_2)& L(p_2,q_3) &L(p_2,q_4)\\
 	L(p_3,q_1)       & L(p_3,q_2)        &L(p_3,q_3) & L(p_3,q_4)\\
 	L(p_4,q_1)& L(p_4,q_2)&L(p_4,q_3)&L(p_4,q_4).\\
 \end{pmatrix}.
\end{equation}
We have
\begin{align}
	\begin{aligned}
	L(p_i,q_j) &= \frac{1}{4\pi}	\int_{[\frac{i-1}{4},\frac{i}{4}]}\int_{[\frac{j-1}{4},\frac{j}{4}]}\frac{\det(\dot{\gamma_1}(s),\dot{\gamma_2}(t),\gamma_1(s)-\gamma_2(t))}{|\gamma_1(s)-\gamma_2(t)|^3}\ ds\ dt\\
	&=\frac{1}{4\pi}	\int_{[\frac{(i-1)\pi}{2},\frac{i\pi}{2}]}\int_{[\frac{(j-1)\pi}{2},\frac{j\pi}{2}]}\frac{\cos(\phi)-\cos(\theta)-\cos(\phi)\cos(\theta)}{(3+2\cos(\phi)-2\cos(\theta)-2\cos(\phi)\cos(\theta))^{\frac{3}{2}}}\ d\phi\ d\theta.\\
\end{aligned}
\end{align}
Numerical calculations show that 
\begin{equation}
	G_1 \approx \begin{pmatrix}
		-0.0640&-0.1413&-0.1413&-0.0640\\
		\quad0.0193&-0.0640&-0.0640&\quad0.0193\\
		\quad0.0193&-0.0640&-0.0640&\quad0.0193\\
		-0.0640&-0.1413&-0.1413&-0.0640\\
		\end{pmatrix}
\end{equation}
and the grand sum of $G_1$ equals   -1.000.
Similarly, we have
\begin{equation}
		G_2 \approx \begin{pmatrix}
		-0.0391&-0.0579&-0.1083&-0.1083&-0.0579&-0.0391\\
		\quad0.0137&\quad0.0069&-0.0653&-0.0653&\quad0.0069&\quad0.0137	\\
			\quad0.0137&\quad0.0069&-0.0653&-0.0653&\quad0.0069&\quad0.0137\\
		-0.0391&-0.0579&-0.1083&-0.1083&-0.0579&-0.0391\\
	\end{pmatrix}
\end{equation}
and the grand sum of $G_2$ equals  -1.000.

\end{example}

\subsection{Multiscale Gauss linking integral}
	\begin{figure}[!htb]
	\centering
	\includegraphics[width=0.5\textwidth]{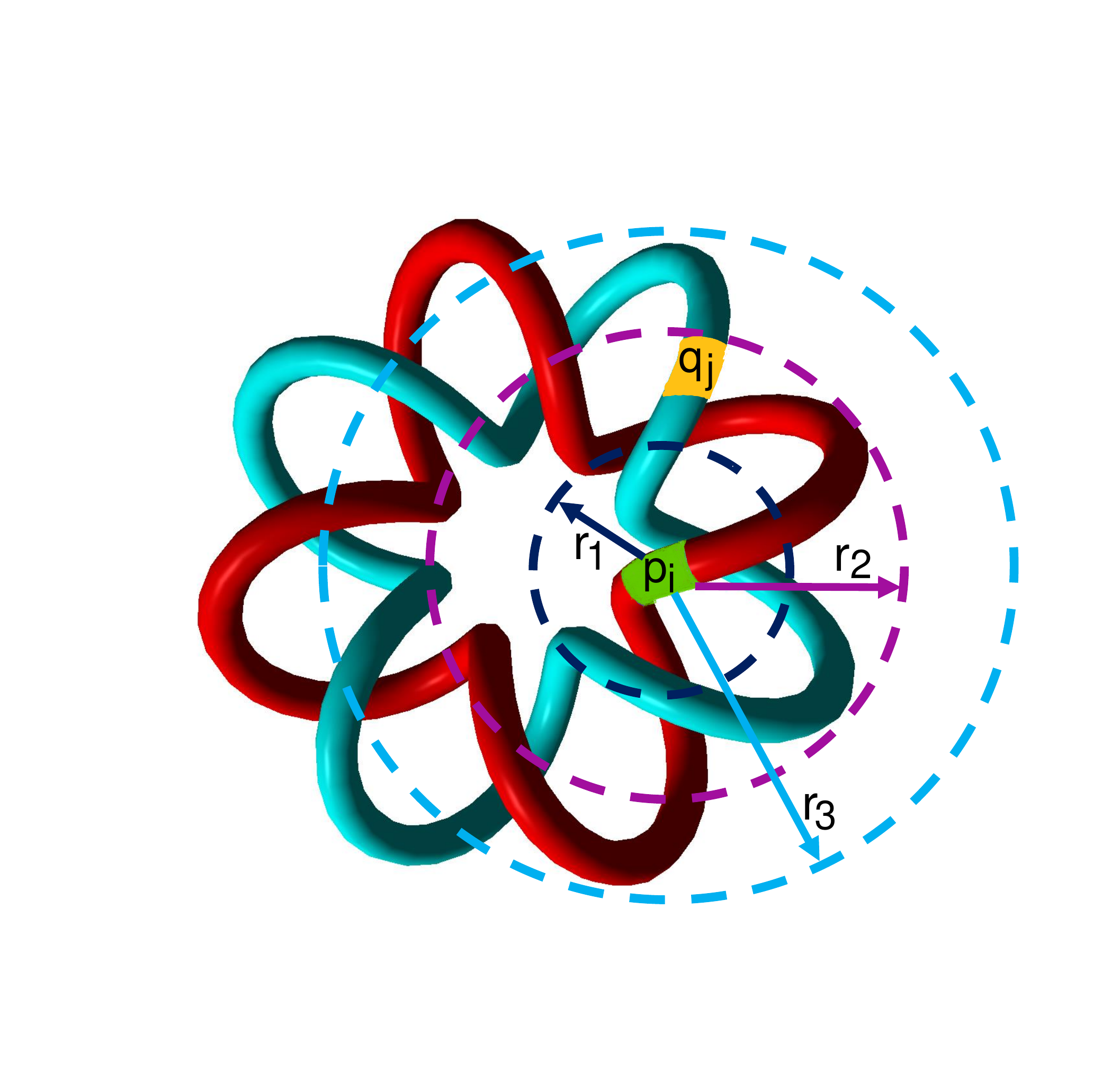}
	\caption{An illustration of multiscale Gauss linking integral with a (2,8) torus link.}
	\label{multiscale}
\end{figure}

In this section, we introduce the multiscale Gauss linking integral (mGLI) based on the segmentation  of the Gauss linking integral. The multiscale Gauss linking integral   gives rise to a vectorization of the Gauss linking integral in the multiscale representation, which can be employed as a  practical embedding for mathematical measurements in various applications. 

\begin{remark} \label{distance}
Multiscale analysis  requires a definition of distance between objects. Since the objects in the  segmentation  of Gauss linking integral are curve segments, we need to define the distance of curve segments $d(p_i,q_j)$. Although distance can be defined in various metrics, we choose the Euclidean distance of two points in this study (see more details in \autoref{eq:dis1} and \autoref{eq:dis2}). For simplicity, here we simply denote the distance as $d(p_i,q_j)$.
\end{remark}

\begin{definition}[Scaled Gauss linking integral]
	Giving a finite set of real numbers $R = \{r_0,r_1,r_2,r_3,\cdots,r_k \}$ where $0=r_0 < r_1<r_2<\cdots<r_k$, the Gauss linking integral at scale $[r_t,r_{t+1}]$ is defined:
	\begin{equation}
			G^{r_t,r_{t+1}} = \begin{pmatrix}
				\chi_{[r_t,r_{t+1}]}(d(p_1,q_1))L(p_1,q_1)&\chi_{[r_t,r_{t+1}]}(d(p_1,q_2))L(p_1,q_2)& \cdots &\chi_{[r_t,r_{t+1}]}(d(p_1,q_m))L(p_1,q_m) \\
				\chi_{[r_t,r_{t+1}]}(d(p_2,q_1))\L(p_2,q_1)&\chi_{[r_t,r_{t+1}]}(d(p_2,q_2))L(p_2,q_2)& \cdots &\chi_{[r_t,r_{t+1}]}(d(p_2,q_m))L(p_2,q_m)\\
				\vdots       & \vdots        &\ddots & \vdots \\
				\chi_{[r_t,r_{t+1}]}(d(p_n,q_1))L(p_n,q_1)& \chi_{[r_t,r_{t+1}]}(d(p_n,q_2))L(p_n,q_2)&\cdots&\chi_{[r_t,r_{t+1}]}(d(p_n,q_m))L(p_n,q_m)\\ 		
			\end{pmatrix},
	\end{equation}
where 
\begin{equation}
	\chi_{[r_t,r_{t+1}]}(x) = \begin{cases}
		1, \text{\rm  if }~ x \in [r_t,r_{t+1}]\\
		0, \text{\rm else}
		\end{cases}
\end{equation}
\end{definition}
\begin{remark}
The scaled Gauss linking integral is used to extract appropriate linking integral within the scale. As shown in \autoref{multiscale}, we have $G_{ij}^{0,r_1}=0$, $G_{ij}^{r_1,r_2}=L(p_i,q_j)$, and $G_{ij}^{r_2,r_3}=0$. By examining multiple Gaussian linking integral at different scales, we obtain an embedding of Gauss linking integral in the multiscale representation. 
\end{remark}

\begin{definition}[Localized scaled   Gauss linking integral]  
For given scale $[r_t,r_{t+1}]$, we can define the localized scaled   Gauss linking integral at $p_i$ or $q_j$ by the followings:
	\begin{equation}
		J^{r_t,r_{t+1}}(p_i) = \sum_{s=1}^{m}G^{r_t,r_{t+1}}_{is}
	\end{equation}
\begin{equation}
	J^{r_t,r_{t+1}}(q_j) = \sum_{s=1}^{n}G^{r_t,r_{t+1}}_{sj}
\end{equation}
\end{definition}
\begin{remark}
The localized scaled Gauss linking integral give rise to a measurement for each curve segment in the curve. By considering different scales, the localized scaled Gauss linking integral provide a featurization of each curve segment $u$:
\begin{equation}\label{feature}
	Feature(u) = (J^{r_1,r_2}(u),J^{r_2,r_3}(u),\cdots,J^{r_{k-1},r_k}(u)).
\end{equation}
The resulting feature vectors offer a basis for machine learning, including deep learning, of complex geometric systems. 
\end{remark}

\subsection{Generalized multiscale Gauss linking integral}

 Vassiliev measure, a generalization of Gauss linking integral,  can be applied to open and closed curves in 3-space \cite{panagiotou2021vassiliev}. Similarly, the proposed mGLI obtained by combining the Gauss linking integral and multiscale process can naturally be applied to links, linkoids, open and closed curves, and other segmentable objects. It can be noticed that any element in the segmentation of the Gauss linking integral is defined on local curve segments. This indicates that one can define a generalized form of the multiscale Gauss linking integral whenever the  segmentation  of the Gauss linking integral is well-defined on local curve segments. In fact, for any topological or geometric structure that can be  segmented  into curve segments $P_n,Q_m$, we can define the following segmentation matrix:

\begin{equation}
	\bar{G} = \begin{pmatrix}
		g(p_1,q_1)&g(p_1,q_2)& \cdots &g(p_1,q_m) \\
		g(p_2,q_1)&g(p_2,q_2)& \cdots &g(p_2,q_m)\\
		\vdots       & \vdots        &\ddots & \vdots \\
		g(p_n,q_1)& g(p_n,q_2)&\cdots&g(p_n,q_m)\\			
	\end{pmatrix},
\end{equation}
where 
\begin{equation}
	g(p_i,q_j)= \begin{cases}
		L(p_i,q_j)& \text{\rm if } p_i \cap q_j \text{ is a null-set,}\\
		0 &\text{ \rm else.}
	\end{cases}
\end{equation}
\begin{figure}[!htb]
	\centering
	\includegraphics[width=0.40\textwidth]{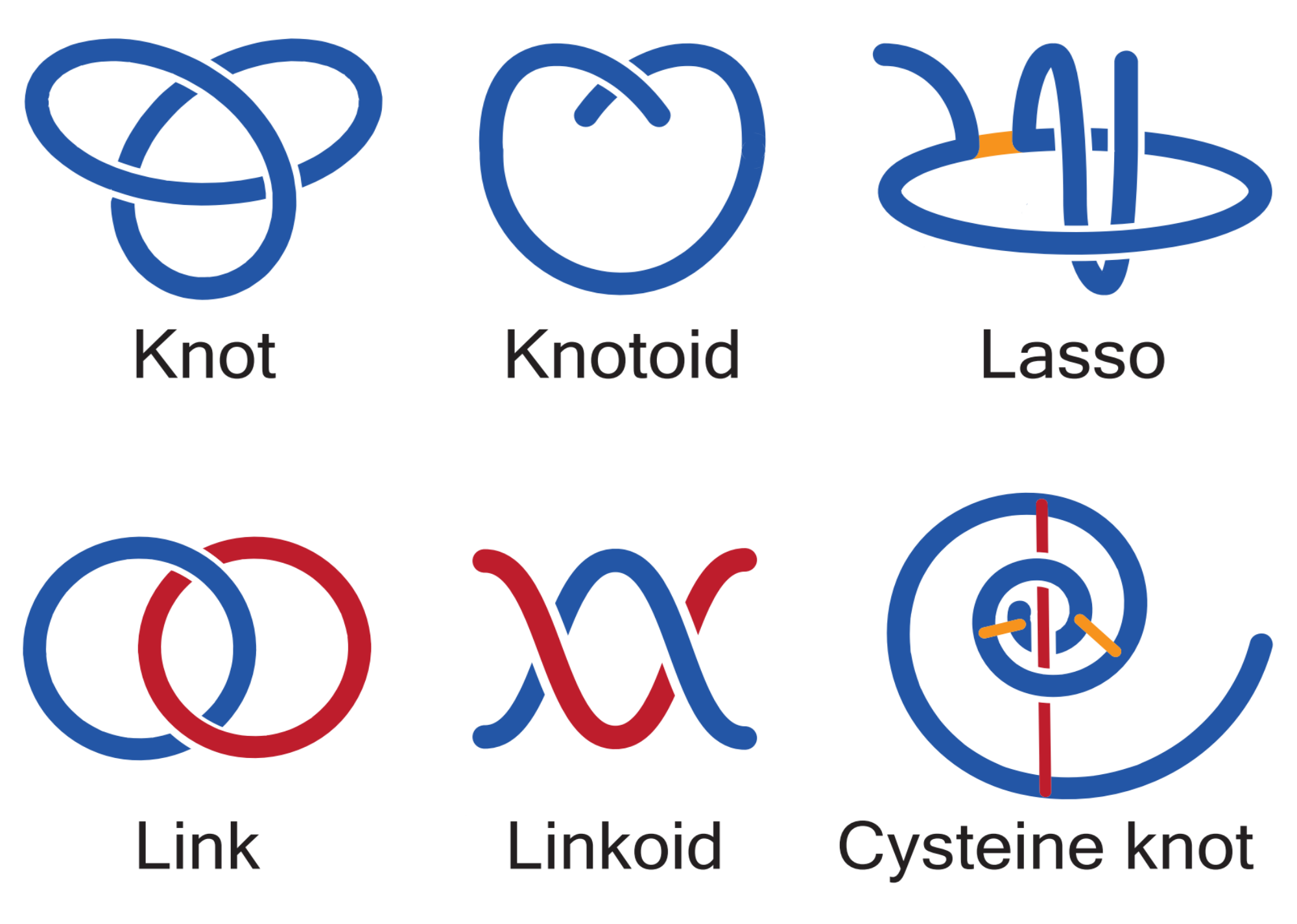}
	\caption{Examples of topological structures which can be studied by multiscale Gauss linking integral. Note that a knotoid is also known as a slipknot \cite{jamroz2015knotprot} in some context.}
	\label{Topo_structure}
\end{figure}
In the above definition, unlike in \autoref{G}, the curve segments in $P_n$ and $Q_m$ are allowed to intersect or even be equal. Thus, the mGLI can be applied in multiple topological/geometric  structures as long as they can locally  be represented as  curve segments. These kinds of structures can be commonly found in multiple studies. Here, we present six examples in \autoref{Topo_structure} including knot, knotoid \cite{gugumcu2017new}, link, linkoid \cite{panagiotou2021vassiliev}, lasso \cite{niemyska2016complex}, and cysteine knot \cite{dabrowski2019knotprot}.

\begin{example}
	\begin{figure}[!htb]
		\centering
		\includegraphics[width=0.15\textwidth]{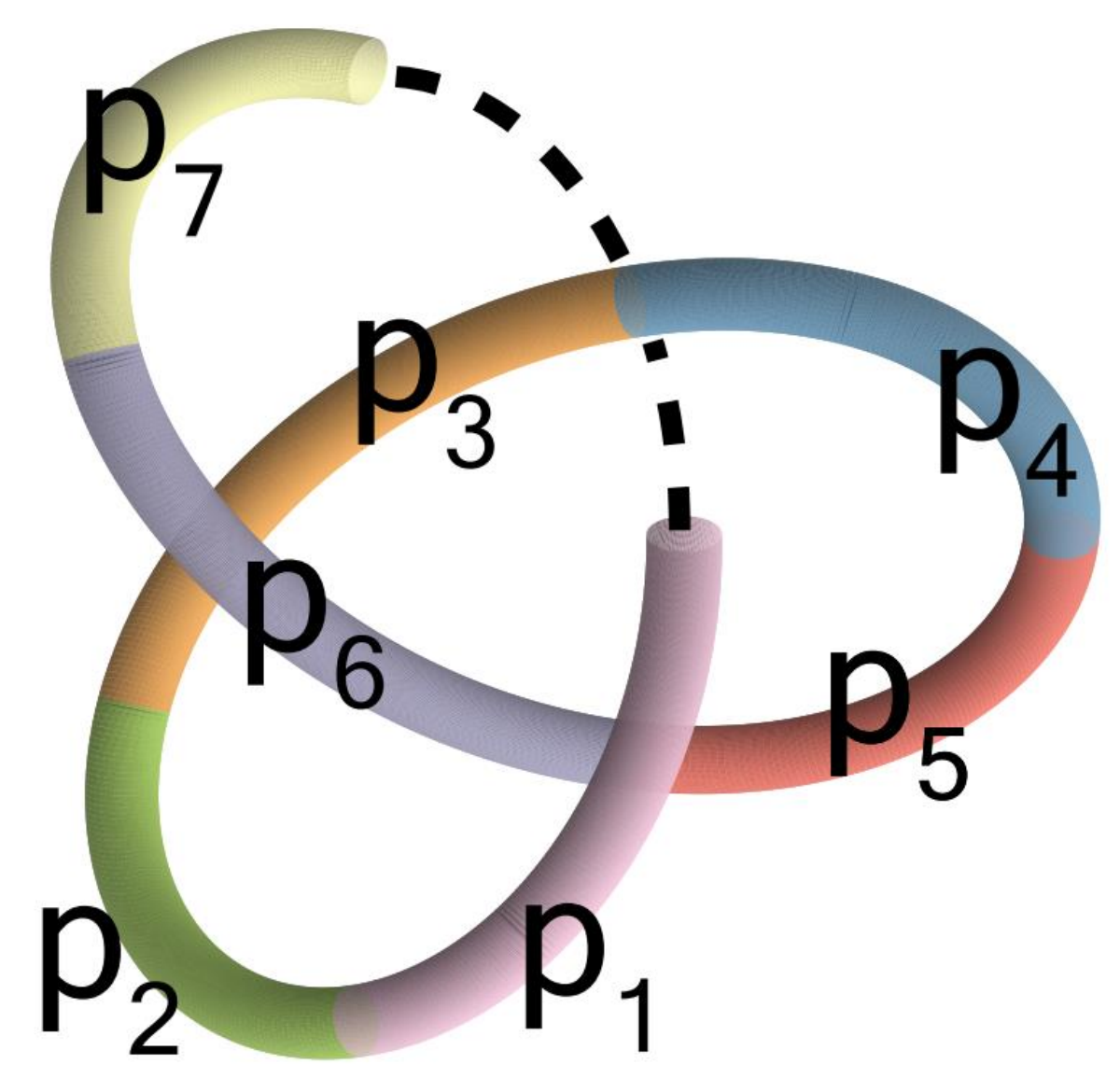}
		\caption{A slipknot with seven curve segments.}
		\label{slipknot}
	\end{figure}
We exhibits in this example how the segmentation  of Gauss linking integral performs on a slipknot shown in \autoref{slipknot}. Here, we consider a curve segmentation $P_7$ with itself, giving the following matrix:
\begin{equation}
	\bar{G}_1 =  \begin{pmatrix}
		g(p_1,p_1)&g(p_1,p_2)& \cdots &g(p_1,p_7) \\
		g(p_2,p_1)&g(p_2,p_2)& \cdots &g(p_2,p_7)\\
		\vdots       & \vdots        &\ddots & \vdots \\
		g(p_7,p_1)& g(p_7,p_2)&\cdots&g(p_7,p_7)\\			
	\end{pmatrix},
\end{equation}
Note that $p_i \cap p_j$ is a null-set if and only if $i \ne j$, and thus,  we have 
\begin{equation}
	g(p_i,p_j)= \begin{cases}
	L(p_i,p_j)& i \ne j\\
	0 & i =j
    \end{cases}.
\end{equation}
Then, we have 
\begin{equation}
	\bar{G}_1=  \begin{pmatrix}
		0 &L(p_1,p_2)& \cdots &L(p_1,p_7) \\
		L(p_2,p_1)&0& \cdots &L(p_2,p_7)\\
		\vdots       & \vdots        &\ddots & \vdots \\
		L(p_7,p_1)& L(p_7,p_2)&\cdots&0\\			
	\end{pmatrix}.
\end{equation}
\end{example}

\begin{example}
		\begin{figure}[!htb]
		\centering
		\includegraphics[width=0.15\textwidth]{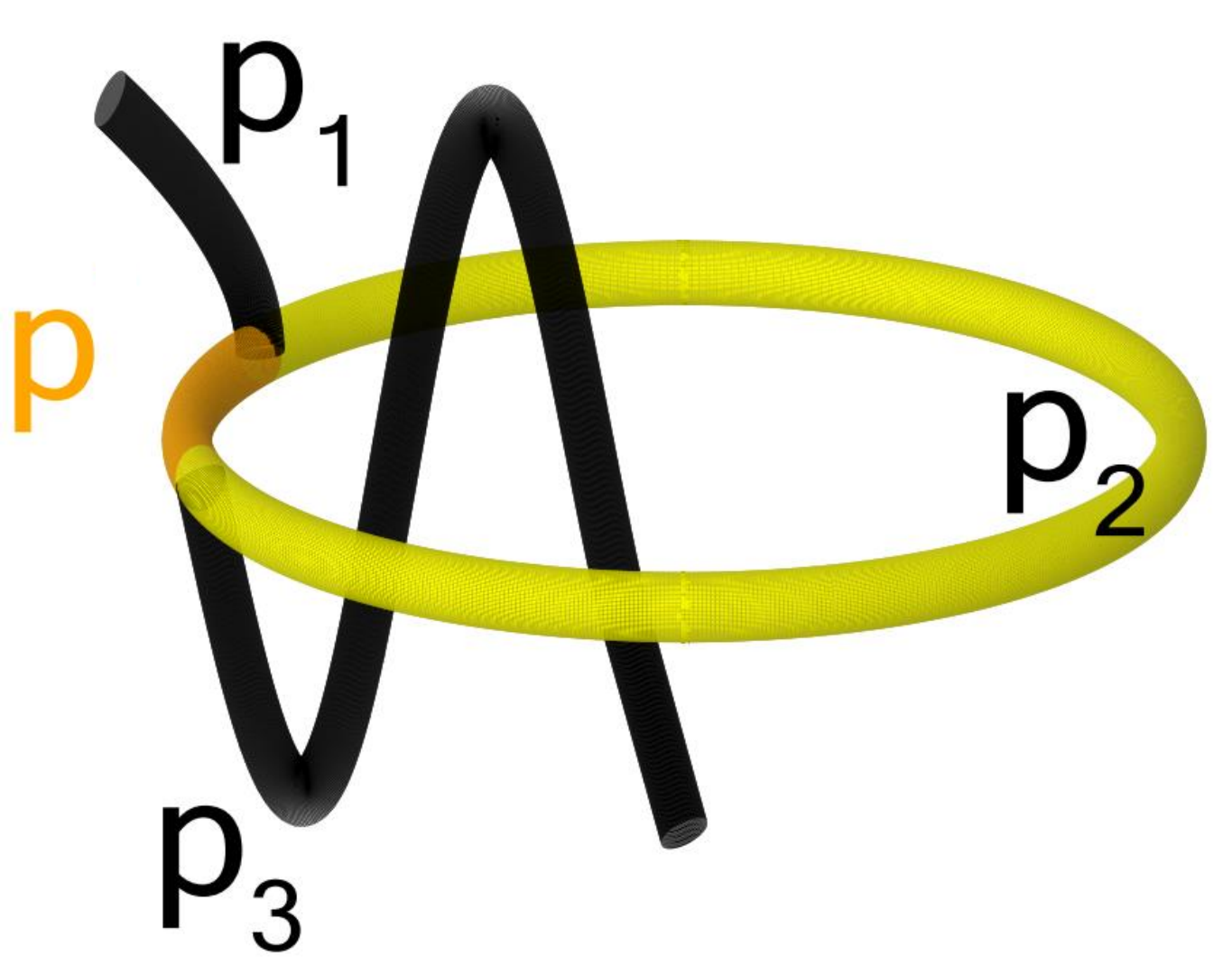}
		\caption{A lasso with four curve segments.}
		\label{lasso}
	\end{figure}
		
	In  this example, we consider a lasso that has a disulfide bridge (orange $p^\prime$) with a slipknot (black $p_1$, yellow $p_2$, and black $p_3$). The   segmentation  of the Gauss linking integral at this lasso is:
	\begin{equation}
			\begin{aligned}
		\bar{G}_2 &=  \begin{pmatrix}
	g(p_1,p_1)&g(p_1,p_2)& g(p_1,p_3) &g(p_1,p^\prime) \\
	g(p_2,p_1)&g(p_2,p_2)& g(p_2,p_3) &g(p_2,p^\prime)\\
	g(p_3,p_1) & g(p_3,p_2) &g(p_3,p_3) & g(p_3,p^\prime) \\
	g(p^\prime,p_1)& g(p^\prime,p_2)&g(p^\prime,p_3)&g(p^\prime,p^\prime)\\			
\end{pmatrix}\\
&=\begin{pmatrix}
	0&g(p_1,p_2)& L(p_1,p_3) &L(p_1,p^\prime) \\
	L(p_2,p_1)&0& L(p_2,p_3) &L(p_2,p^\prime)\\
	L(p_3,p_1) & L(p_3,p_2) &0 & L(p_3,p^\prime) \\
	L(p^\prime,p_1)& L(p^\prime,p_2)&L(p^\prime,p_3)&0\\			
\end{pmatrix}.
\end{aligned}
	\end{equation}
\end{example}

\section{Numerical experiment}
\subsection{Protein segmentation and multiscale analysis}
\begin{figure}[!htb]
		\includegraphics[width=\textwidth]{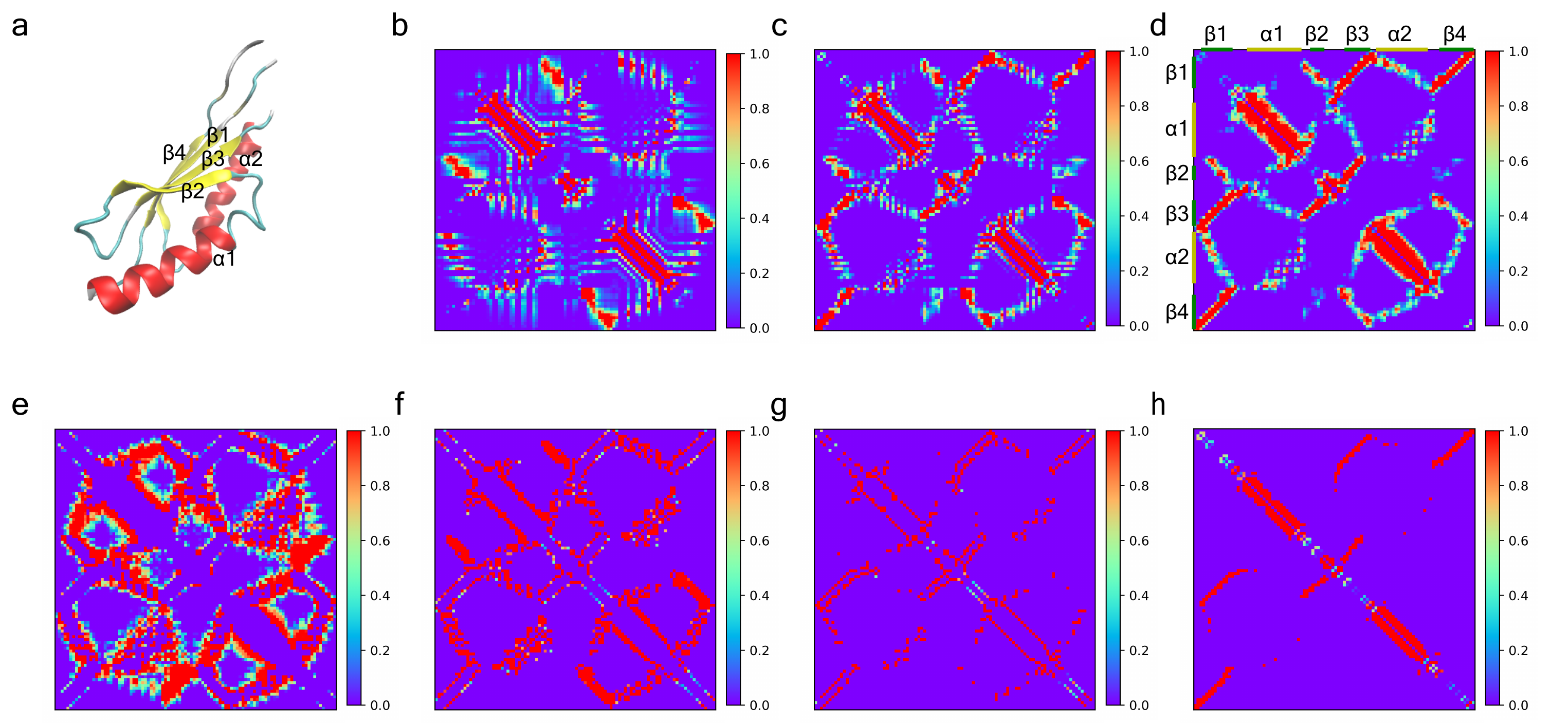}
		\caption{An illustration of multiscale Guass linking integral on a protein (PDBID: 1J27). {\bf a} The 3D structure of 1J27  consisting of two $\alpha$-helices and four $\beta$-sheets. {\bf b} The segmentation of the Gauss linking integral of protein 1J27. {\bf c} The absolute value of Gauss linking integral matrix of protein 1J27. {\bf d} The absolute Gauss linking integral matrix of protein 1J27. {\bf e}-{\bf h} Absolute Gauss linking integral matrices of protein 1J27 at different scales.}
		\label{concept}
	\end{figure}
In this section, we assess the effectiveness of the proposed generalized mGLI in predicting protein B-factors. Protein B-factors, also known as temperature factors or atomic displacement parameters, provide information about the mobility and flexibility of atoms in a protein structure. B-factors are derived from crystallographic experiments or other experimental means and are represented as a numerical value associated with each atom in the protein. High B-factors indicate increased atomic mobility, suggesting regions of the protein that are flexible or undergoing conformational changes. Low B-factors, on the other hand, indicate rigid regions with limited atomic motion.

To predict protein B-factors, we consider a protein as an open curve, acknowledging that the polypeptide chain of a protein molecule can be seen as an open polygon $l$ whose vertices are corresponding to the C$_\alpha$ atoms, while the edges represent the pseudobonds that connect a C$_\alpha$  atom to another one in an adjacent amino acid residue. 

Here, we propose a curve segmentation induced by C$_\alpha$ atoms:
\begin{equation} \label{eq:dis1} 
	p_i= \{x \in l_1| f(x,c_i) = \inf_{c\in C}f(x,c)\}, 1\le i \le n,
\end{equation}
where $f(a,b)$ is the distance of points $a$ and $b$ alone $l$, and $C$ is the set of  C$_\alpha$ atoms.
Then, the $d(p_i,q_j)$ assumed in \autoref{distance} can be defined:
\begin{equation}\label{eq:dis2}
	d(p_i,q_j) = d_E(c_i,c_j),
\end{equation}
where $d_E$ is the Euclidean distance in the 3D space. 

 Then, the segmentation of Gauss linking integral  that investigates the inter-crossings between segments of the protein can be given: 
\begin{equation}
	G = \begin{pmatrix}
		\hat{L}(p_1,p_1)&\hat{L}(p_1,p_2)& \cdots &\hat{L}(p_1,p_n) \\
		\hat{L}(p_2,p_1)&\hat{L}(p_2,p_2)& \cdots &\hat{L}(p_2,p_n)\\
		\vdots       & \vdots        &\ddots & \vdots \\
		\hat{L}(p_n,p_1)& \hat{L}(p_n,p_2)&\cdots&\hat{L}(p_n,p_n)\\			
	\end{pmatrix}
     =\begin{pmatrix}
     	0&L(p_1,p_2)& \cdots &L(p_1,p_n) \\
     	L(p_2,p_1)&0& \cdots &L(p_2,p_n)\\
     	\vdots       & \vdots        &\ddots & \vdots \\
     	L(p_n,p_1)& L(p_n,p_2)&\cdots&0\\			
     \end{pmatrix}.
\end{equation}
Here, we present the illustration of Gauss linking integral matrix for a protein (PDBID: 1J27) in \autoref{concept}b. This matrix can reflect much  structural information about the protein. Therefore, the Gauss linking integral matrix itself is an excellent structural featurization method for proteins. 

Also, it is worth noting that the value of Gauss linking integral is subject to its consideration of the curve orientations. On a specific task that does not depend on the orientation of the curve, removing the influence of the orientation may give a better results. In particular, in protein B-factor analysis, the effects of two differently oriented structures on a certain  C$_\alpha$ atom will hardly cancel each other out. In fact, by taking the absolute value of segmentation of Gauss linking integral, the orientation information can be greatly reduced. As shown in \autoref{concept}{\bf c}, a large number of interlaced patterns are removed. Further, to completely ignore orientation information, we consider the absolute Gauss linking integral
\begin{equation}
	\bar{L}(l_1,l_2) = \frac{1}{4\pi}	\int_{[0,1]}\int_{[0,1]}\left|{\frac{\det(\dot{\gamma_1}(s),\dot{\gamma_2}(t),\gamma_1(s)-\gamma_2(t))}{|\gamma_1(s)-\gamma_2(t)|^3}}\right|\ ds\ dt,
\end{equation}
with its corresponding matrix. An illustration of this approach is given in  \autoref{concept}{\bf d}.

To present featurization of the absolute Gauss linking integral, we consider the localized scaled absolute Gauss linking integral at each curve segment around each C$_\alpha$ atom of various scales, generating a  multiscale Gauss linking integral feature vector. We validate in numerous protein structures that there exists a high correlation coefficients between the reciprocal of featurization of segments (\autoref{feature}) and the corresponding C$_\alpha$ atom's B-factor. 

\subsection{Protein flexibility analysis}
	
	Protein flexibility analysis is the study of the dynamic behavior and conformational changes of proteins. Proteins are not rigid structures but rather exhibit flexibility and undergo various motions and deformations essential for their biological functions.	Protein flexibility analysis involves investigating the conformational space accessible to a protein, understanding its inherent flexibility, and exploring how it relates to protein function, stability, and interactions with other molecules. It provides insights into how proteins can adopt different shapes and conformations to perform their roles, such as enzyme catalysis, molecular recognition, and signal transduction.
	
	Several computational methods have been proposed for protein flexiblity analysis. The Gauss network model (GNM) \cite{rader2005gaussian} and anisotropic network model (ANM) \cite{eyal2006anisotropic}, and normal mode analysis (NMA)\cite{bahar2005coarse} have been widely used for their simplicity. However, Park et al. \cite{park2013coarse} have demonstrated that both GNM and NMA were ineffective in analyzing numerous protein structures. Their findings revealed that, on average, the correlation coefficient of GNM and NMA, for the three protein sets, categorized as small-size, medium-size, and large-size, was below 0.6 and 0.5, respectively. Recently, many advanced methods have been proposed for this problem, including flexibility rigidity index-based methods, pfFRI \cite{opron2014fast} and opFRI \cite{opron2014fast} and  topology-based methods, including atom-specific persistent homology (ASPH) \cite{bramer2020atom} and evolutionary homology (EH) \cite{cang2020evolutionary}.
	
\begin{figure}[htb!]
	\centering
	\begin{tabular}{c}
				\includegraphics[width=0.6\textwidth]{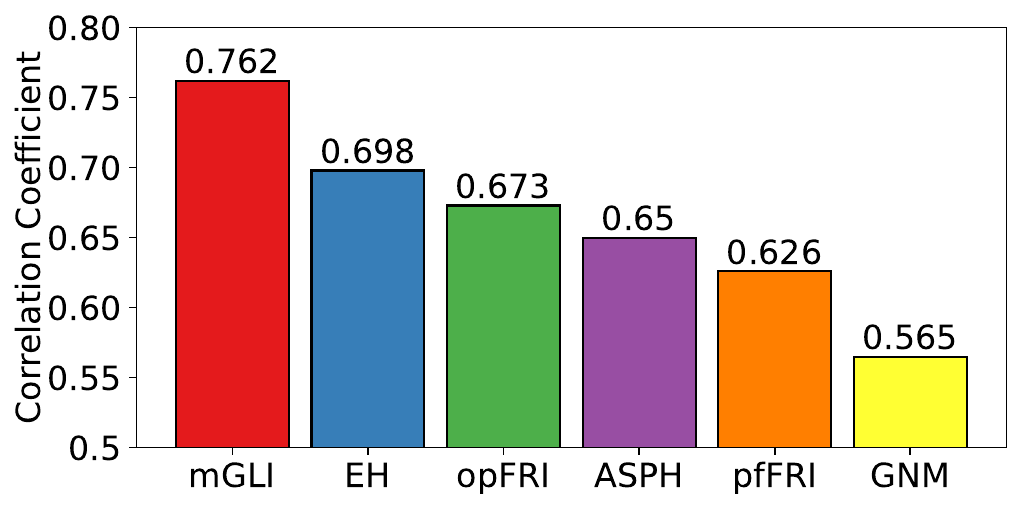}
	\end{tabular}
		\caption{Comparison of B-factor predictions by six methods for 364 proteins. }
\label{Comparison}
\end{figure}	
	
To assess the performance of the proposed multiscale Gauss linking integral (mGLI) for protein flexibility analysis, we utilized a benchmark dataset comprising 364 protein structures sourced from \cite{opron2014fast}. This dataset was used for comparisons against previous methods, namely opFRI \cite{opron2014fast}, pfFRI \cite{opron2014fast}, and GNM \cite{park2013coarse}. In our evaluation, we set the mGLI with a maximal cutoff of 27Å.

In \autoref{Comparison}, we present the comparative results of mGLI with previous methods for each protein in the dataset (refer to \autoref{364} for detailed information). Remarkably, mGLI outperformed the previous methods in 320 out of 364 proteins. On average, mGLI achieved the highest correlation coefficient of 0.762, surpassing the values of 0.674 for opFRI, 0.628 for pfFRI, and 0.567 for GNM. This represents an improvement of 13.1\%, 21.3\%, and 34.4\% respectively. These results strongly suggest that mGLI is a promising method for protein flexibility analysis, as the higher correlation coefficients obtained indicate its effectiveness.
	
\begin{table}[htb!]
	\centering
	\begin{tabular}{|c|ccccccc|}
		\hline
		Protein set &EH\cite{cang2020evolutionary}&ASPH\cite{bramer2020atom}&opFRI\cite{opron2014fast} & pfFRI\cite{opron2014fast} & GNM\cite{park2013coarse} &NMA\cite{park2013coarse} &mGLI\\\hline
		Small&0.773&0.870&0.667&0.594&0.541&0.480& {\bf 0.933}\\
		Medium&0.729&0.680&0.664&0.605&0.550&0.482& {\bf 0.813} \\
		Large&0.665&0.610&0.636&0.591&0.529&0.494& {\bf 0.750} \\ \hline
	\end{tabular}
\caption{Average correlation coefficients for  C$_\alpha$ B-factor prediction of our mGLI method with previous EH, opFRI, pfFRI, GNM, and NMA for three protein sets of different size. }
	\label{comparison_table}
\end{table}

In addition, to validate the effectiveness of mGLI for predicting C$_\alpha$ atom B-factors in proteins of different sizes, we compared our method with previous approaches including EH \cite{cang2020evolutionary}, ASPH \cite{bramer2020atom}, opFRI \cite{opron2014fast}, pfFRI \cite{opron2014fast}, GNM \cite{park2013coarse}, and NMA \cite{park2013coarse} on three protein sets, as shown in \autoref{comparison_table}.

mGLI achieved average correlation coefficients of 0.933, 0.813, and 0.750 for the small, medium, and large protein sets, respectively, in predicting C$_\alpha$ atom B-factors. Our results on the three datasets significantly outperformed the previous methods, demonstrating improvements of 27.3\%, 11.5\%, and 12.7\% on the small, medium, and large protein sets, respectively, compared to the previous state-of-the-art method EH.

The improvement of our method in predicting protein B-factors is primarily attributed to the clever utilization of knot theory. Traditional B-factor analysis methods primarily focus on individual atoms and their spatial positions in three-dimensional space, considering the thermal motion and disorder of atoms within a protein structure. However, the consideration of bonding interactions between atoms, which indirectly influences the observed B-factor values, is rarely utilized in B-factor analysis. By incorporating knot theory, our method introduces the concept of pseudobonds between protein atoms, thereby capturing the influence of bonding interactions.

The combination of knot theory with the multiscale procedure allows for the localization of measurements provided by the knot theory, taking advantage of the spatial positions and environments of atoms. This synergy between multiscale analysis and knot theory results in a powerful method for protein B-factor prediction, highlighting the potential of multiscale approaches in localizing measurements derived from knot theory.

\begin{figure}[!htb]
	\includegraphics[width=0.7\textwidth]{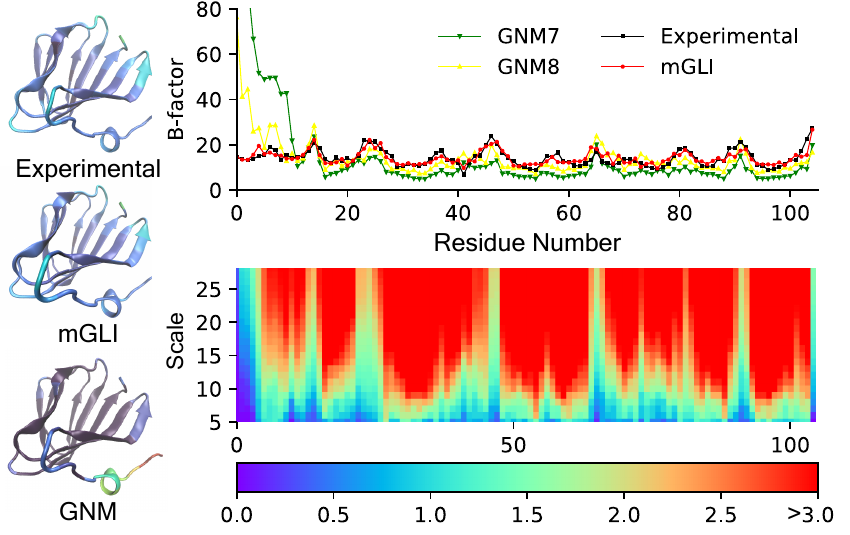}
	\caption{The visualization of mGLI method on the probable antibiotics synthesis protein 1V70. \textbf{Left} are the 3D structures of 1V70, colored by experimental B-factor, mGLI pedicted B-factor, and GNM predicted B-factor, from top to bottom, respetively. \textbf{Right top} is the comparison of GNM7, GNM8, and mGLI prediction values and experimental value for each residue. The $x$-axis represents the residue number, and the $y$-axis represents the B-factor. \textbf{Right bottom} is the vitalization of mGLI features with the maximal cutoff at 30\AA. The $x$-axis represents the residue number and the $y$-axis represents the supreme of filtration range. Note that all value exceed 3.0 are labeled as red. GNM7 represents the GNM with cutoff distance 7\AA.}
	\label{1v70}
\end{figure}

To provide a direct visualization of the multiscale Gauss linking integral (mGLI), we consider an example using a probable antibiotics synthesis protein (PDBID: 1V70) with a residue number of 105, as shown in \autoref{1v70}. Traditional GNM methods exhibit a large error around residues 1-10. In contrast, mGLI successfully predicts the B-factor of this region. The predictions of mGLI are directly obtained from the multiscale Gauss linking integral features.
A visualization of the multiscale Gauss linking integral features has been displayed in the bottom right of \autoref{1v70}. For each scale, we consider the accumulated abosoluted Gauss linking integral. A bar corresponding to the color and the accumulated value is placed below it. Note that we label all value exceed a specific value as red (3.0 in this case). Thus we could notice that the upper bound of scale required to accumulate the absolute Gauss linking integral value of each C$_\alpha$ atom to 3.0 is different. In general, the smaller the upper limit of scale needed to reach a specific value , the smaller the corresponding B-factor value. In fact, the peaks in the blue and green areas at the bottom often coincide with the peaks of the experimental B-factor.

\begin{figure}[!htb] 	
	\includegraphics[width=0.7\textwidth]{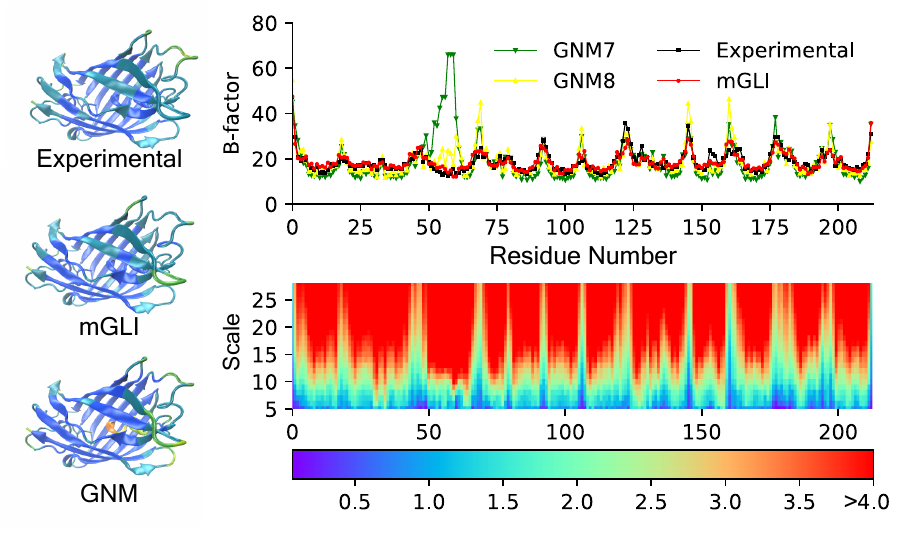}
	\caption{The visualization of mGLI method on the monomeric cyan fluorescent protein 2HQK. \textbf{Left} are the 3D structures of 2HQK, colored by experimental B-factor, mGLI pedicted B-factor, and GNM predicted B-factor, from top to bottom, respectively. \textbf{Right top} is the comparison of GNM7, GNM8, and mGLI prediction value and experimental value for each residue. The $x$-axis represents the residue number, and the $y$-axis represents the B-factor. \textbf{Right bottom} is the virtualization of mGLI features with the maximal cutoff at 30\AA. The $x$-axis represents the residue number and the $y$-axis represents the supreme of scales. Note that all values exceed 4.0 are labeled as red. GNM7 represents the GNM with cutoff distance 7\AA. }
	\label{2hqk}		
\end{figure}

The above observations are quite general. We can further validate these findings by providing another visualization using a monomeric cyan fluorescent protein (PDBID: 2HQK), as shown in \autoref{2hqk}. In this visualization, we label all absolute Gauss linking integral values exceed 4.0 as red. The protein 2HQK has a residue number of 219, which is more than twice that of 1V70. As expected, there are more peaks in its B-factors. Interestingly, each peak in the B-factor plot corresponds nicely to a peak in the visualization of the absolute Gauss linking integral feature. This correlation reinforces the relationship between the B-factor values and the absolute Gauss linking integral. Furthermore, the absolute Gauss linking integral feature demonstrates how mGLI avoids the erroneous predictions made by traditional GNM methods in the area around residues 50-60. This showcases the advantage of mGLI in accurately capturing the protein's flexibility and providing more reliable predictions.
By examining these different protein examples, we consistently observe the correlation between the absolute Gauss linking integral feature and the experimental B-factor values, highlighting the effectiveness of mGLI in protein B-factor prediction.

\section{Discussions}
	
\subsection{Absolute multiscale Gauss linking integral}
	
	\label{sec_dis}
\begin{figure}[!htb]
		\includegraphics[width=0.6\textwidth]{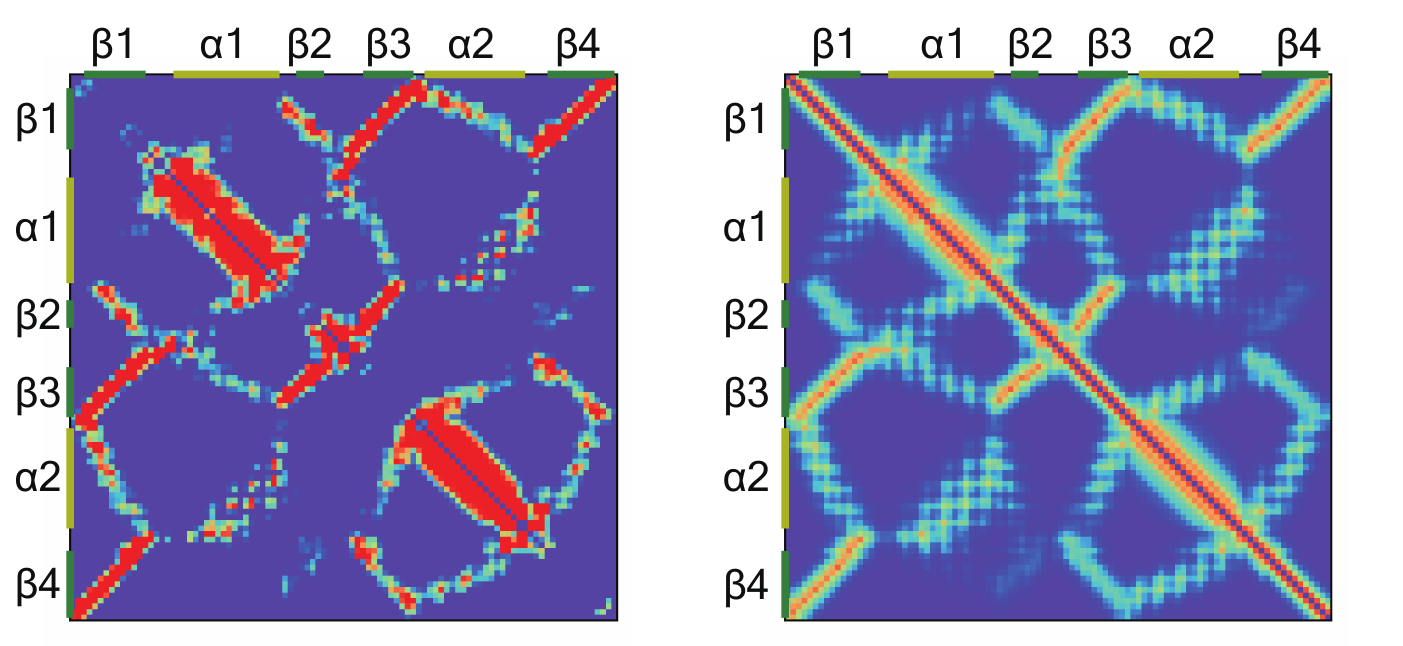}
		\label{dis_compare}
		\caption{\textbf{Left} is the absolute Gauss linking integral matrix for protein 1J27. \textbf{Right} is the transient probability matrix (TPM) for protein 1J27 \cite{opron2015communication}. Each point in top and left are colored green or yellow, if the point represents a residue that  is in a $\beta$-sheet or an $\alpha$-helix of 1J27.}
\end{figure}

The Gauss linking integral provides information that previous methods cannot capture, which is one of the key factors contributing to the success of mGLI. To illustrate the information contained in mGLI, we compare the segmentation of the absolute multiscale Gauss linking integral with the previous transient probability matrix (TPM) \cite{opron2015communication} for protein 1J27.
The structural information that can be observed through the previous method becomes more evident and clear in mGLI. For example, in the TPM, the $\alpha 1$-$\alpha 1$ and $\alpha 2$-$\alpha 2$ interactions are represented as slightly thicker yellow blocks along the diagonal. In contrast, in mGLI, these interactions are depicted as more expressive and prominent large red blocks. This enhanced visualization allows for a better distinction of the self-interaction of the alpha chain from other structural elements, such as the self-interaction of the region between $\beta2$ and $\beta3$.
Moreover, the contrast between different values within each block is stronger in mGLI compared to TPM. This is particularly evident in the blocks representing interactions such as $\alpha 1$-$\alpha 2$, $\beta 1$-$\alpha 2$, $\beta 1$-$\beta 2$, and so on. In the multiscale analysis, the differences in mGLI values at different scales are further amplified, providing insights into varying structural information.
Overall, mGLI offers a more expressive and informative representation of the protein structure compared to the previous TPM method, enabling a better understanding of the intricate interactions within the protein.

\subsection{Multiscale featurization}

\begin{figure}[!htb]
	\centering
	\includegraphics[width=0.25 \textwidth]{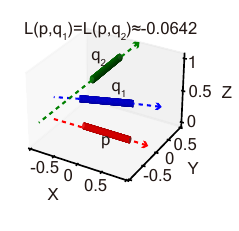}
	\caption{An example of line segments in different scales reporting the same Gauss linking integral.}
	\label{example2}
\end{figure}
The process of multiscale vectorizing Gauss linking integral in the multiscale analysis allows mGLI to extract a wealth of information directly related to the spatial arrangement of segments. Each scale captures the Gauss linking integral at a different level of granularity, enabling the discrimination of various segment arrangements.

For example, \autoref{example2} illustrates two different segment arrangements of $q$ in 3-dimensional space that have identical Gauss linking integrals with $p$ (i.e., $L(p,q_1)=L(p,q_2)$). Without the multiscale process, it would be challenging to differentiate between these arrangements based on the Gauss linking integral alone. However, by considering different scales such as $[0,0.5]$ and $[0.5,1]$, the arrangements can be distinguished.

\subsection{Curve segmentation and multiscale granularity}

In principle, our method allows for the arbitrary combination of curve segmentations with any multiscale schemes. However, in practical applications, the performance of mGLI is highly dependent on the selection of curve segmentations and multiscale schemes.
First and foremost, the values of the Gauss linking integral of a local curve segment depend not only on their spatial alignment but also on their relative lengths compared to the global curve. When the length of a curve segment approaches zero, the corresponding Gauss linking integral approximates to 0. Similarly, as the curve segments expand to cover the global curve, the Gauss linking integral returns global information. In both cases, the Gauss linking integral fails to extract useful spatial information regarding local alignments.
Secondly, the selection of the multiscale scheme impacts the featurization of the Gauss linking integral. Ideally, different scales should capture different spatial structure information. If the information difference between scales is too small, it can result in a large number of identical or trivial features. Conversely, if the difference is too large, it may lead to the loss of details.

In the specific context of our protein flexibility analysis, we naturally choose a segment that precisely covers a single $C_\alpha$ atom along the polymer chain. Additionally, in our study, the multiscale scheme is selected to start from  5\AA~and extend up to 27\AA, with each scale interval set at 1\AA. This choice is based on the fact that the average distance between $C_\alpha$ atoms is approximately 3.8\AA. Such a selection of the multiscale scheme results in a powerful featurization method that provides abundant representations of local protein structures.

\section{Conclusion}	

In recent years, topological data analysis (TDA) has experienced rapid growth across various scientific and engineering disciplines. The central technique of TDA, persistent homology, utilizes a multiscale filtration process to enrich algebraic topology. However, the impact of knot theory and related mathematical approaches on data science has been limited, partly due to their global formalism.

This work introduces knot data analysis (KDA) as a powerful new paradigm for data science. Our approach involves partitioning knots and other curved objects into segments and conducting a multiscale analysis at each segment. At each scale, multiscale Gauss link integrals (mGLIs) are defined to capture the local structure and connectivity. These mGLIs allow for the measurement of the properties of knots and links at different scales, ultimately recovering their global properties at a sufficiently large scale. Additionally, we extend mGLIs to encompass various segmentable objects, such as knots, knotoids, lassos, links, linkoids, cysteine knots, and more.

While persistent homology remains a global technique despite filtration, our proposed KDA offers a local approach that focuses on delineating the localized environment around each segment at different scales. By adopting this localized perspective, KDA provides valuable insights into the specific characteristics and relationships within the data.

To showcase the effectiveness of the proposed KDA, we apply the multiscale Gauss link integral (mGLI) to protein flexibility prediction. By utilizing mGLI, we can quantify the relationship between C$_\alpha$ atoms at different scales, enabling a comprehensive analysis of protein flexibility.

In our study, we conducted experiments using four benchmark datasets and compared the performance of mGLI with several previous methods, including a persistent homology-based approach. The results of these extensive experiments consistently demonstrate the superiority of mGLI. Specifically, mGLI achieves a minimum improvement of 13\% compared to the previous state-of-the-art methods. Notably, it outperforms the persistent homology-based method by approximately 21\%.

These results highlight the significant potential of the proposed KDA in the field of data science. By leveraging the power of mGLI and its ability to capture multiscale relationships, KDA offers promising opportunities for improving various data analysis tasks.

Since mGLI offers multiscale representations of a wide range of spatial structure and connectivity information, it can be leveraged by machine learning algorithms, enabling artificial intelligence to extract powerful and effective insights into segmentable structures.
The rich information captured by mGLI allows machine learning algorithms to explore and understand the intricate relationships and patterns within the data. This enables AI systems to make more accurate predictions, classify complex structures, and uncover hidden features that may not be apparent through traditional approaches.


\section*{Acknowledgment}	
			This work was supported in part by NIH grants R01GM126189, R01AI164266, and R35GM148196
, NSF grants
			DMS-2052983, DMS-1761320, DMS-2245903,  and IIS-1900473, NASA grant 80NSSC21M0023, MSU Foundation, Bristol-Myers Squibb 65109, and Pfizer.

	\newpage
 
\section{Appendix}

In this appendix, detailed B-factor predictions for a set of 364 proteins are presented. 
Comparison is given to  opFRI\cite{opron2014fast},  pfFRI\cite{opron2014fast},    GNM\cite{park2013coarse},    and mGLI.

\begin{longtable}{|l|r|r|r|r|r|l|r|r|r|r|r|}

	\hline
	PDBID &    N &  opFRI\cite{opron2014fast} &  pfFRI\cite{opron2014fast} & GNM\cite{park2013coarse} &  mGLI & PDBID &    N &  opFRI &  pfFRI &    GNM &  mGLI \\ \hline
	\hline
	\endfirsthead
	
	\hline
	PDBID &    N &  opFRI &  pfFRI &    GNM &  mGLI & PDBID &    N &  opFRI &  pfFRI &    GNM &  mGLI \\ \hline
	\hline
	\endhead
	\endfoot

  1ABA &   87 &          0.727 &          0.698 &          0.613 & \textbf{0.856} &  1PEF &   18 &          0.888 &          0.826 &          0.808 &   \textbf{1.0} \\ \hline
 1AHO &   64 &          0.698 &          0.625 &          0.562 & \textbf{0.788} &  1PEN &   16 &          0.516 &          0.465 &           0.27 & \textbf{0.858} \\ \hline
 1AIE &   31 &          0.588 &          0.416 &          0.155 & \textbf{0.985} &  1PMY &  123 &          0.671 &          0.654 &          0.685 & \textbf{0.763} \\ \hline
 1AKG &   16 &          0.373 &           0.35 &          0.185 & \textbf{0.854} &  1PZ4 &  114 &          0.828 &          0.781 &          0.843 & \textbf{0.925} \\ \hline
 1ATG &  231 &          0.613 &          0.578 &          0.497 & \textbf{0.743} &  1Q9B &   43 &          0.746 &          0.726 &          0.656 & \textbf{0.985} \\ \hline
 1BGF &  124 &          0.603 &          0.539 &          0.543 & \textbf{0.921} &  1QAU &  112 &          0.678 &          0.672 &           0.62 & \textbf{0.731} \\ \hline
 1BX7 &   51 &          0.726 &          0.623 &          0.706 & \textbf{0.893} &  1QKI & 3912 & \textbf{0.809} &          0.751 &          0.645 &          0.678 \\ \hline
 1BYI &  224 &          0.543 &          0.491 &          0.552 & \textbf{0.678} &  1QTO &  122 &          0.543 &           0.52 &          0.334 & \textbf{0.725} \\ \hline
 1CCR &  111 &           0.58 &          0.512 &          0.351 & \textbf{0.643} &  1R29 &  122 &           0.65 &          0.631 &          0.556 & \textbf{0.781} \\ \hline
 1CYO &   88 &          0.751 &          0.702 &          0.741 &  \textbf{0.91} &  1R7J &   90 &          0.789 &          0.621 &          0.368 & \textbf{0.929} \\ \hline
 1DF4 &   57 &          0.912 &          0.889 &          0.832 & \textbf{0.975} &  1RJU &   36 &          0.517 &          0.447 &          0.431 & \textbf{0.839} \\ \hline
 1E5K &  188 &          0.746 &          0.732 & \textbf{0.859} &          0.831 &  1RRO &  112 &          0.435 &          0.372 &          0.529 & \textbf{0.691} \\ \hline
 1ES5 &  260 &          0.653 &          0.638 &          0.677 & \textbf{0.681} &  1SAU &  114 &          0.742 &          0.671 &          0.596 & \textbf{0.837} \\ \hline
 1ETL &   12 &           0.71 &          0.609 &          0.628 &   \textbf{1.0} &  1TGR &  104 &           0.72 &          0.711 &          0.714 & \textbf{0.805} \\ \hline
 1ETM &   12 &          0.544 &          0.393 &          0.432 &  \textbf{0.97} &  1TZV &  141 &          0.837 &           0.82 &          0.841 & \textbf{0.915} \\ \hline
 1ETN &   12 &          0.089 &          0.023 &         -0.274 &  \textbf{0.67} &  1U06 &   55 &          0.474 &          0.429 &          0.434 & \textbf{0.893} \\ \hline
 1EW4 &  106 &           0.65 &          0.644 &          0.547 & \textbf{0.797} &  1U7I &  267 & \textbf{0.778} &          0.762 &          0.691 &          0.758 \\ \hline
 1F8R & 1932 & \textbf{0.878} &          0.859 &          0.738 &          0.746 &  1U9C &  221 &            0.6 &          0.577 &          0.522 & \textbf{0.754} \\ \hline
 1FF4 &   65 &          0.718 &          0.613 &          0.674 & \textbf{0.871} &  1UHA &   83 &          0.726 &          0.665 &          0.638 & \textbf{0.863} \\ \hline
 1FK5 &   93 &           0.59 &          0.568 &          0.485 & \textbf{0.779} &  1UKU &  102 &          0.665 &          0.661 &          0.742 & \textbf{0.847} \\ \hline
 1GCO & 1044 & \textbf{0.766} &          0.693 &          0.646 &          0.706 &  1ULR &   87 &          0.639 &          0.594 &          0.495 & \textbf{0.822} \\ \hline
 1GK7 &   39 &          0.845 &          0.773 &          0.821 & \textbf{0.935} &  1UOY &   64 &          0.713 &          0.653 &          0.671 & \textbf{0.877} \\ \hline
 1GVD &   52 &          0.781 &          0.732 &          0.591 & \textbf{0.933} &  1USE &   40 &          0.438 &          0.146 &         -0.142 & \textbf{0.964} \\ \hline
 1GXU &   88 &          0.748 &          0.634 &          0.421 & \textbf{0.882} &  1USM &   77 &          0.832 &          0.809 &          0.798 & \textbf{0.964} \\ \hline
 1H6V & 2927 &          0.488 &          0.429 &          0.306 & \textbf{0.504} &  1UTG &   70 &          0.691 &           0.61 &          0.538 & \textbf{0.849} \\ \hline
 1HJE &   13 &          0.811 &          0.686 &          0.616 & \textbf{0.975} &  1V05 &   96 &          0.629 &          0.599 &          0.632 & \textbf{0.742} \\ \hline
 1I71 &   83 &          0.549 &          0.516 &          0.549 & \textbf{0.691} &  1V70 &  105 &          0.622 &          0.492 &          0.162 & \textbf{0.794} \\ \hline
 1IDP &  441 & \textbf{0.735} &          0.715 &           0.69 &          0.729 &  1VRZ &   21 &          0.792 &          0.695 &          0.677 & \textbf{0.999} \\ \hline
 1IFR &  113 &          0.697 &          0.689 &          0.637 &  \textbf{0.77} &  1W2L &   97 &          0.691 &          0.564 &          0.397 & \textbf{0.845} \\ \hline
 1K8U &   89 &          0.553 &          0.531 &          0.378 & \textbf{0.763} &  1WBE &  204 &          0.591 &          0.577 &          0.549 & \textbf{0.706} \\ \hline
 1KMM & 1499 & \textbf{0.749} &          0.744 &          0.558 &          0.727 &  1WHI &  122 & \textbf{0.601} &          0.539 &           0.27 &          0.591 \\ \hline
 1KNG &  144 &          0.547 &          0.536 &          0.512 & \textbf{0.723} &  1WLY &  322 &          0.695 &          0.679 &          0.666 & \textbf{0.701} \\ \hline
 1KR4 &  110 &          0.635 &          0.612 &          0.466 & \textbf{0.693} &  1WPA &  107 &          0.634 &          0.577 &          0.417 & \textbf{0.793} \\ \hline
 1KYC &   15 &          0.796 &          0.763 &          0.754 &   \textbf{1.0} &  1X3O &   80 &            0.6 &          0.559 &          0.654 &  \textbf{0.72} \\ \hline
 1LR7 &   73 &          0.679 &          0.657 &           0.62 & \textbf{0.882} &  1XY1 &   18 &          0.832 &          0.645 &          0.447 &   \textbf{1.0} \\ \hline
 1MF7 &  194 &          0.687 &          0.681 &            0.7 & \textbf{0.717} &  1XY2 &    8 &          0.619 &           0.57 &          0.562 &   \textbf{1.0} \\ \hline
 1N7E &   95 &          0.651 &          0.609 &          0.497 & \textbf{0.794} &  1Y6X &   87 &          0.596 &          0.524 &          0.366 & \textbf{0.739} \\ \hline
 1NKD &   59 &           0.75 &          0.703 &          0.631 & \textbf{0.913} &  1YJO &    6 &          0.375 &          0.333 &          0.434 &   \textbf{1.0} \\ \hline
 1NKO &  122 &          0.619 &          0.535 &          0.368 & \textbf{0.703} &  1YZM &   46 &          0.842 &          0.834 &          0.901 & \textbf{0.975} \\ \hline
 1NLS &  238 &          0.669 &           0.53 &          0.523 & \textbf{0.777} &  1Z21 &   96 &          0.662 &          0.638 &          0.433 & \textbf{0.739} \\ \hline
 1NNX &   93 &          0.795 &          0.789 &          0.631 & \textbf{0.903} &  1ZCE &  146 &          0.808 &          0.757 &           0.77 & \textbf{0.911} \\ \hline
 1NOA &  113 &          0.622 &          0.604 &          0.615 & \textbf{0.629} &  1ZVA &   75 &          0.756 &          0.579 &           0.69 & \textbf{0.934} \\ \hline
 1NOT &   13 &          0.746 &          0.622 &          0.523 & \textbf{0.999} &  2A50 &  457 & \textbf{0.564} &          0.524 &          0.281 &          0.527 \\ \hline
 1O06 &   20 &           0.91 &          0.874 &          0.844 &   \textbf{1.0} &  2AGK &  233 &          0.705 &          0.694 &          0.512 & \textbf{0.753} \\ \hline
 1O08 &  221 &          0.562 &          0.333 &          0.309 &   \textbf{0.8} &  2AH1 &  939 & \textbf{0.684} &          0.593 &          0.521 &          0.557 \\ \hline
 1OB4 &   16 &          0.776 &          0.763 &           0.75 &   \textbf{1.0} &  2B0A &  186 &          0.639 &          0.603 &          0.467 & \textbf{0.664} \\ \hline
 1OB7 &   16 &          0.737 &          0.545 &          0.652 &   \textbf{1.0} &  2BCM &  413 & \textbf{0.555} &          0.551 &          0.477 &          0.442 \\ \hline
 1OPD &   85 &          0.555 &          0.409 &          0.398 & \textbf{0.689} &  2BF9 &   36 &          0.606 &          0.554 &           0.68 & \textbf{0.999} \\ \hline
 1P9I &   29 &          0.754 &          0.742 &          0.625 & \textbf{0.975} &  2BRF &  100 &          0.795 &          0.764 &           0.71 & \textbf{0.857} \\ \hline
 2CE0 &   99 &          0.706 &          0.598 &          0.529 &  \textbf{0.85} &  2C71 &  205 &          0.658 &          0.649 &           0.56 & \textbf{0.848} \\ \hline
 2CG7 &   90 &          0.551 &          0.539 &          0.379 & \textbf{0.743} &  2OLX &    4 &          0.917 &          0.888 &          0.885 &   \textbf{1.0} \\ \hline
 2COV &  534 & \textbf{0.846} &          0.823 &          0.812 &          0.638 &  2PKT &   93 &          0.162 &          0.003 &          0.193 & \textbf{0.641} \\ \hline
 2CWS &  227 &          0.647 &           0.64 &          0.696 & \textbf{0.773} &  2PLT &   99 &          0.508 &          0.484 &          0.509 & \textbf{0.663} \\ \hline
 2D5W & 1214 & \textbf{0.689} &          0.682 &          0.681 &          0.648 &  2PMR &   76 &          0.693 &          0.682 &          0.619 & \textbf{0.878} \\ \hline
 2DKO &  253 & \textbf{0.816} &          0.812 &           0.69 &          0.789 &  2POF &  440 & \textbf{0.682} &          0.651 &          0.589 &          0.637 \\ \hline
 2DPL &  565 &          0.596 &          0.538 & \textbf{0.658} &          0.479 &  2PPN &  107 &          0.677 &          0.638 &          0.668 & \textbf{0.765} \\ \hline
 2DSX &   52 &          0.337 &          0.333 &          0.127 &  \textbf{0.65} &  2PSF &  608 &          0.526 &            0.5 & \textbf{0.565} &          0.552 \\ \hline
 2E10 &  439 & \textbf{0.798} &          0.796 &          0.692 &          0.777 &  2PTH &  193 &          0.822 &          0.784 &          0.767 & \textbf{0.837} \\ \hline
 2E3H &   81 &          0.692 &          0.682 &          0.605 & \textbf{0.722} &  2Q4N &  153 &          0.711 &          0.667 &           0.74 &  \textbf{0.79} \\ \hline
 2EAQ &   89 &          0.753 &           0.69 &          0.695 & \textbf{0.862} &  2Q52 &  412 & \textbf{0.756} &          0.748 &          0.621 &           0.71 \\ \hline
 2EHP &  248 & \textbf{0.804} & \textbf{0.804} &          0.773 &          0.762 &  2QJL &   99 &          0.594 &          0.584 &          0.594 & \textbf{0.812} \\ \hline
 2EHS &   75 &           0.72 &          0.713 &          0.747 & \textbf{0.851} &  2R16 &  176 &          0.582 &          0.495 &          0.618 & \textbf{0.675} \\ \hline
 2ERW &   53 &          0.461 &          0.253 &          0.199 & \textbf{0.858} &  2R6Q &  138 &          0.603 &           0.54 &          0.529 & \textbf{0.653} \\ \hline
 2ETX &  389 &           0.58 &          0.556 & \textbf{0.632} &          0.621 &  2RB8 &   93 &          0.727 &          0.614 &          0.517 & \textbf{0.808} \\ \hline
 2FB6 &  116 &          0.791 &          0.786 &           0.74 & \textbf{0.838} &  2RE2 &  238 &          0.652 &          0.613 &          0.673 & \textbf{0.683} \\ \hline
 2FG1 &  157 &           0.62 &          0.617 &          0.584 & \textbf{0.741} &  2RFR &  154 &          0.693 &          0.671 &          0.753 & \textbf{0.816} \\ \hline
 2FN9 &  560 &          0.607 &          0.595 &          0.611 & \textbf{0.659} &  2V9V &  135 &          0.555 &          0.548 &          0.528 & \textbf{0.679} \\ \hline
 2FQ3 &   85 &          0.719 &          0.692 &          0.348 & \textbf{0.868} &  2VE8 &  515 & \textbf{0.744} &          0.643 &          0.616 &          0.569 \\ \hline
 2G69 &   99 &          0.622 &           0.59 &          0.436 & \textbf{0.812} &  2VH7 &   94 &          0.775 &          0.726 &          0.596 & \textbf{0.897} \\ \hline
 2G7O &   68 &          0.785 &          0.784 &           0.66 & \textbf{0.952} &  2VIM &  104 &          0.413 &          0.393 &          0.212 & \textbf{0.598} \\ \hline
 2G7S &  190 &           0.67 &          0.644 &          0.649 &  \textbf{0.76} &  2VPA &  204 & \textbf{0.763} &          0.755 &          0.576 &          0.641 \\ \hline
 2GKG &  122 &          0.688 &          0.646 &          0.711 & \textbf{0.859} &  2VQ4 &  106 &           0.68 &          0.679 &          0.555 & \textbf{0.861} \\ \hline
 2GOM &  121 & \textbf{0.586} &          0.584 &          0.491 &          0.461 &  2VY8 &  149 &  \textbf{0.77} &          0.724 &          0.533 &          0.701 \\ \hline
 2GXG &  140 &          0.847 &           0.78 &           0.52 &   \textbf{0.9} &  2VYO &  210 &          0.675 &          0.648 &          0.729 & \textbf{0.767} \\ \hline
 2GZQ &  191 &          0.505 &          0.382 &          0.369 &  \textbf{0.86} &  2W1V &  548 &  \textbf{0.68} &  \textbf{0.68} &          0.571 &           0.66 \\ \hline
 2HQK &  213 &          0.824 &          0.809 &          0.365 &  \textbf{0.89} &  2W2A &  350 & \textbf{0.706} &          0.638 &          0.589 &          0.631 \\ \hline
 2HYK &  238 &          0.585 &          0.575 &           0.51 & \textbf{0.626} &  2W6A &  117 & \textbf{0.823} &          0.748 &          0.647 &          0.761 \\ \hline
 2I24 &  113 &          0.593 &          0.498 &          0.494 & \textbf{0.711} &  2WJ5 &   96 &          0.484 &           0.44 &          0.357 & \textbf{0.656} \\ \hline
 2I49 &  398 &          0.714 &          0.683 &          0.601 & \textbf{0.739} &  2WUJ &  100 & \textbf{0.739} &          0.598 &          0.598 &          0.668 \\ \hline
 2IBL &  108 &          0.629 &          0.625 &          0.352 & \textbf{0.723} &  2WW7 &  150 & \textbf{0.499} &          0.471 &          0.356 &           0.46 \\ \hline
 2IGD &   61 &          0.585 &          0.481 &          0.386 & \textbf{0.825} &  2WWE &  111 &          0.692 &          0.582 &          0.628 & \textbf{0.807} \\ \hline
 2IMF &  203 &          0.652 &          0.625 &          0.514 & \textbf{0.719} &  2X1Q &  240 &          0.534 &          0.478 &          0.443 & \textbf{0.549} \\ \hline
 2IP6 &   87 &          0.654 &          0.578 &          0.572 & \textbf{0.862} &  2X25 &  168 &          0.632 &          0.598 &          0.403 & \textbf{0.793} \\ \hline
 2IVY &   88 &          0.544 &          0.483 &          0.271 & \textbf{0.838} &  2X3M &  166 &          0.744 &          0.717 &          0.655 & \textbf{0.832} \\ \hline
 2J32 &  244 & \textbf{0.863} &          0.848 &          0.855 &          0.802 &  2X5Y &  171 &          0.718 &          0.705 &          0.694 & \textbf{0.811} \\ \hline
 2J9W &  200 &          0.716 &          0.705 &          0.662 & \textbf{0.745} &  2X9Z &  262 &          0.583 &          0.578 &          0.574 &  \textbf{0.64} \\ \hline
 2JKU &   35 &          0.805 &          0.695 &          0.656 & \textbf{0.993} &  2XHF &  310 &          0.606 &          0.591 &          0.569 & \textbf{0.721} \\ \hline
 2JLI &  100 &          0.779 &          0.613 &          0.622 & \textbf{0.836} &  2Y0T &  101 &          0.778 &          0.774 &          0.798 & \textbf{0.835} \\ \hline
 2JLJ &  115 &          0.741 &           0.72 &          0.527 & \textbf{0.844} &  2Y72 &  170 &           0.78 &          0.754 &          0.766 & \textbf{0.874} \\ \hline
 2MCM &  113 &          0.789 &          0.713 &          0.639 & \textbf{0.846} &  2Y7L &  319 & \textbf{0.928} &          0.797 &          0.747 &          0.742 \\ \hline
 2NLS &   36 &          0.605 &          0.559 &           0.53 & \textbf{0.942} &  2Y9F &  149 &          0.771 &          0.762 &          0.664 & \textbf{0.835} \\ \hline
 2NR7 &  194 &          0.803 &          0.785 &          0.727 & \textbf{0.868} &  2YLB &  400 & \textbf{0.807} & \textbf{0.807} &          0.675 &            0.7 \\ \hline
 2NUH &  104 &          0.835 &          0.691 &          0.771 & \textbf{0.908} &  2YNY &  315 & \textbf{0.813} &          0.804 &          0.706 &          0.427 \\ \hline
 2O6X &  306 & \textbf{0.814} &          0.799 &          0.651 &          0.797 &  2ZCM &  357 & \textbf{0.458} &          0.422 &           0.42 &          0.434 \\ \hline
 2OA2 &  132 &          0.571 &          0.456 &          0.458 & \textbf{0.794} &  2ZU1 &  360 & \textbf{0.689} &          0.672 &          0.653 &          0.649 \\ \hline
 2OCT &  192 &          0.567 &           0.55 &           0.54 & \textbf{0.675} &  3A0M &  148 & \textbf{0.807} &          0.712 &          0.392 &          0.704 \\ \hline
 2OHW &  256 &          0.614 &          0.539 &          0.475 &  \textbf{0.74} &  3A7L &  128 &          0.713 &          0.663 & \textbf{0.756} &          0.733 \\ \hline
 2OKT &  342 &          0.433 &          0.411 &          0.336 & \textbf{0.517} &  3AMC &  614 &          0.675 &          0.669 &          0.581 & \textbf{0.726} \\ \hline
 2OL9 &    6 &          0.909 &          0.904 &          0.689 &   \textbf{1.0} &  3AUB &  116 &          0.614 &          0.608 &          0.637 & \textbf{0.674} \\ \hline
 3BA1 &  312 &          0.661 &          0.624 &          0.621 & \textbf{0.687} &  3B5O &  230 & \textbf{0.644} &          0.629 &          0.601 &          0.614 \\ \hline
 3BED &  261 & \textbf{0.845} &           0.82 &          0.684 &          0.769 &  3MD4 &   12 &           0.86 &          0.781 &          0.914 &   \textbf{1.0} \\ \hline
 3BQX &  139 &          0.634 &          0.481 &          0.297 &  \textbf{0.72} &  3MD5 &   12 &          0.649 &          0.413 &         -0.218 &   \textbf{1.0} \\ \hline
 3BZQ &   99 &          0.532 &          0.516 &          0.466 & \textbf{0.717} &  3MEA &  166 &          0.669 &          0.669 &            0.6 & \textbf{0.694} \\ \hline
 3BZZ &  100 &          0.485 &           0.45 &            0.6 & \textbf{0.616} &  3MGN &  348 &          0.205 &          0.119 &          0.193 & \textbf{0.267} \\ \hline
 3DRF &  547 &          0.559 &          0.549 &          0.488 & \textbf{0.601} &  3MRE &  383 & \textbf{0.661} &          0.641 &          0.567 &          0.513 \\ \hline
 3DWV &  325 & \textbf{0.707} &          0.661 &          0.547 &          0.702 &  3N11 &  325 & \textbf{0.614} &          0.583 &          0.517 &          0.608 \\ \hline
 3E5T &  228 &          0.502 &          0.489 &          0.296 & \textbf{0.601} &  3NE0 &  208 &          0.706 &          0.645 &          0.659 & \textbf{0.811} \\ \hline
 3E7R &   40 &          0.706 &          0.687 &          0.642 & \textbf{0.929} &  3NGG &   94 &          0.696 &          0.689 &          0.719 & \textbf{0.847} \\ \hline
 3EUR &  140 &          0.431 &          0.427 & \textbf{0.577} &           0.52 &  3NPV &  495 & \textbf{0.702} &          0.653 &          0.677 &          0.574 \\ \hline
 3F2Z &  149 &          0.824 &          0.792 &           0.74 & \textbf{0.861} &  3NVG &    6 &          0.721 &          0.617 &          0.597 &   \textbf{1.0} \\ \hline
 3F7E &  254 & \textbf{0.812} &          0.803 &          0.811 &          0.809 &  3NZL &   73 &          0.627 &          0.583 &          0.506 & \textbf{0.891} \\ \hline
 3FCN &  158 &           0.64 &          0.606 &          0.632 & \textbf{0.847} &  3O0P &  194 &          0.727 &          0.706 &          0.734 & \textbf{0.772} \\ \hline
 3FE7 &   91 &          0.583 &          0.533 &          0.276 & \textbf{0.685} &  3O5P &  128 &          0.734 &          0.698 &           0.63 & \textbf{0.902} \\ \hline
 3FKE &  250 &          0.525 &          0.476 &          0.435 & \textbf{0.769} &  3OBQ &  150 &          0.649 &          0.645 & \textbf{0.655} &          0.647 \\ \hline
 3FMY &   66 &          0.701 &          0.655 &          0.556 & \textbf{0.876} &  3OQY &  234 &          0.698 &          0.686 &          0.637 & \textbf{0.725} \\ \hline
 3FOD &   48 &          0.532 &           0.44 &         -0.126 & \textbf{0.768} &  3P6J &  125 &          0.774 &          0.767 &           0.81 & \textbf{0.893} \\ \hline
 3FSO &  221 & \textbf{0.831} &          0.817 &          0.793 &            0.6 &  3PD7 &  188 &           0.77 &          0.723 &          0.589 & \textbf{0.778} \\ \hline
 3FTD &  240 &          0.722 &          0.713 &          0.634 & \textbf{0.733} &  3PES &  165 & \textbf{0.697} &          0.642 &          0.683 &          0.686 \\ \hline
 3FVA &    6 &          0.835 &          0.825 &          0.789 &   \textbf{1.0} &  3PID &  387 &          0.537 &          0.531 & \textbf{0.642} &          0.575 \\ \hline
 3G1S &  418 & \textbf{0.771} &            0.7 &           0.63 &          0.653 &  3PIW &  154 &          0.758 &          0.744 &          0.717 & \textbf{0.861} \\ \hline
 3GBW &  161 &           0.82 &          0.747 &           0.51 & \textbf{0.839} &  3PKV &  221 &          0.625 &          0.597 &          0.568 & \textbf{0.673} \\ \hline
 3GHJ &  116 &          0.732 &          0.511 &          0.196 & \textbf{0.813} &  3PSM &   94 & \textbf{0.876} &           0.79 &          0.745 &           0.74 \\ \hline
 3HFO &  197 & \textbf{0.691} &           0.67 &          0.518 &          0.646 &  3PTL &  289 & \textbf{0.543} &          0.541 &          0.468 &          0.535 \\ \hline
 3HHP & 1234 &  \textbf{0.72} &          0.716 &          0.683 &          0.643 &  3PVE &  347 & \textbf{0.718} &          0.667 &          0.568 &          0.581 \\ \hline
 3HNY &  156 &          0.793 &          0.723 &          0.758 &  \textbf{0.89} &  3PZ9 &  357 &          0.709 &          0.709 &          0.678 & \textbf{0.803} \\ \hline
 3HP4 &  183 &          0.534 &            0.5 &          0.573 & \textbf{0.677} &  3PZZ &   12 &          0.945 &          0.922 &           0.95 &   \textbf{1.0} \\ \hline
 3HWU &  144 &          0.754 &          0.748 &          0.841 & \textbf{0.861} &  3Q2X &    6 &          0.922 &          0.904 &          0.866 &   \textbf{1.0} \\ \hline
 3HYD &    7 &          0.966 &           0.95 &          0.867 &   \textbf{1.0} &  3Q6L &  131 &          0.622 &          0.577 &          0.605 & \textbf{0.699} \\ \hline
 3HZ8 &  192 &          0.617 &          0.502 &          0.475 & \textbf{0.634} &  3QDS &  284 &  \textbf{0.78} &          0.745 &          0.568 &          0.742 \\ \hline
 3I2V &  124 &          0.486 &          0.441 &          0.301 & \textbf{0.562} &  3QPA &  197 &          0.587 &          0.442 &          0.503 & \textbf{0.866} \\ \hline
 3I2Z &  138 & \textbf{0.613} &          0.599 &          0.317 &          0.597 &  3R6D &  221 &          0.688 &          0.669 &          0.495 & \textbf{0.774} \\ \hline
 3I4O &  135 &          0.735 &          0.714 & \textbf{0.738} &          0.732 &  3R87 &  132 &          0.452 &          0.419 &          0.286 & \textbf{0.701} \\ \hline
 3I7M &  134 &          0.667 &          0.635 &          0.695 & \textbf{0.774} &  3RQ9 &  162 &  \textbf{0.51} &          0.403 &          0.242 &          0.507 \\ \hline
 3IHS &  169 & \textbf{0.586} &          0.565 &          0.409 &          0.583 &  3RY0 &  128 &          0.616 &          0.606 &           0.47 &  \textbf{0.77} \\ \hline
 3IVV &  149 & \textbf{0.817} &          0.797 &          0.693 &          0.734 &  3RZY &  139 &            0.8 &          0.784 &          0.849 & \textbf{0.879} \\ \hline
 3K6Y &  227 &          0.586 &          0.535 &          0.301 & \textbf{0.629} &  3S0A &  119 &          0.562 &          0.524 &          0.526 & \textbf{0.672} \\ \hline
 3KBE &  140 &          0.705 &          0.704 &          0.611 & \textbf{0.726} &  3SD2 &   86 &          0.523 &          0.421 &          0.237 & \textbf{0.734} \\ \hline
 3KGK &  190 & \textbf{0.784} &          0.775 &           0.68 &          0.662 &  3SEB &  238 &          0.801 &          0.712 &          0.826 & \textbf{0.871} \\ \hline
 3KZD &   85 &          0.647 &          0.611 &          0.475 & \textbf{0.719} &  3SED &  124 &          0.709 &          0.658 &          0.712 &  \textbf{0.82} \\ \hline
 3L41 &  220 & \textbf{0.718} &          0.716 &          0.669 &          0.678 &  3SO6 &  150 &          0.675 &          0.666 &           0.63 & \textbf{0.848} \\ \hline
 3LAA &  169 & \textbf{0.827} &          0.647 &          0.659 &           0.82 &  3SR3 &  637 &          0.619 &          0.611 &          0.624 &  \textbf{0.67} \\ \hline
 3LAX &  106 &          0.734 &           0.73 &          0.584 & \textbf{0.898} &  3SUK &  248 & \textbf{0.644} &          0.633 &          0.567 &          0.619 \\ \hline
 3LG3 &  833 & \textbf{0.658} &          0.614 &          0.589 &          0.652 &  3SZH &  697 & \textbf{0.817} &          0.815 &          0.697 &          0.706 \\ \hline
 3LJI &  272 &          0.612 &          0.608 &          0.551 & \textbf{0.619} &  3T0H &  208 &          0.808 &          0.775 &          0.694 & \textbf{0.889} \\ \hline
 3M3P &  249 &          0.584 &          0.554 &          0.338 & \textbf{0.757} &  3T3K &  122 &          0.796 &          0.748 &          0.735 & \textbf{0.918} \\ \hline
 3M8J &  178 &  \textbf{0.73} &          0.728 &          0.628 &          0.674 &  3T47 &  141 &          0.592 &          0.527 &          0.447 & \textbf{0.649} \\ \hline
 3M9J &  210 & \textbf{0.639} &          0.574 &          0.296 &          0.499 &  3TDN &  357 &          0.458 &          0.419 &           0.24 & \textbf{0.584} \\ \hline
 3M9Q &  176 &          0.591 &           0.51 &          0.471 & \textbf{0.627} &  3TOW &  152 &          0.578 &          0.556 &          0.571 & \textbf{0.802} \\ \hline
 3MAB &  173 & \textbf{0.664} &          0.591 &          0.451 &          0.589 &  3TUA &  210 &          0.665 &          0.658 &          0.588 & \textbf{0.815} \\ \hline
 3U6G &  248 & \textbf{0.635} &          0.632 &          0.526 &          0.527 &  3TYS &   75 &          0.853 &            0.8 &          0.791 & \textbf{0.939} \\ \hline
 3U97 &   77 &          0.753 &          0.736 &          0.712 & \textbf{0.916} &  4DT4 &  160 &          0.776 &          0.738 &          0.716 & \textbf{0.813} \\ \hline
 3UCI &   72 &          0.589 &          0.526 &          0.495 & \textbf{0.844} &  4EK3 &  287 &  \textbf{0.68} &  \textbf{0.68} &          0.674 &          0.676 \\ \hline
 3UR8 &  637 & \textbf{0.666} &          0.652 &          0.597 &          0.651 &  4ERY &  318 &           0.74 &          0.701 &          0.688 &  \textbf{0.78} \\ \hline
 3US6 &  148 &          0.698 &          0.586 &          0.553 & \textbf{0.878} &  4ES1 &   95 &          0.648 &          0.625 &          0.551 & \textbf{0.813} \\ \hline
 3V1A &   48 &          0.531 &          0.487 &          0.583 & \textbf{0.978} &  4EUG &  225 &           0.57 &          0.529 &          0.405 & \textbf{0.617} \\ \hline
 3V75 &  285 &          0.604 &          0.596 &          0.491 &  \textbf{0.71} &  4F01 &  448 &          0.633 &          0.372 & \textbf{0.688} &          0.377 \\ \hline
 3VN0 &  193 &  \textbf{0.84} &          0.837 &          0.812 &          0.734 &  4F3J &  143 &          0.617 &          0.598 &          0.551 &   \textbf{0.7} \\ \hline
 3VOR &  182 &          0.602 &          0.557 &          0.484 & \textbf{0.674} &  4FR9 &  141 &          0.671 &          0.655 &          0.501 & \textbf{0.807} \\ \hline
 3VUB &  101 &          0.625 &           0.61 &          0.607 &  \textbf{0.84} &  4G14 &   15 &          0.467 &          0.323 &          0.356 &   \textbf{1.0} \\ \hline
 3VVV &  108 &          0.833 &          0.741 &          0.753 & \textbf{0.899} &  4G2E &  151 &           0.76 &          0.755 &          0.758 & \textbf{0.863} \\ \hline
 3VZ9 &  163 & \textbf{0.785} &          0.749 &          0.695 &          0.725 &  4G5X &  550 & \textbf{0.786} &          0.754 &          0.743 &          0.775 \\ \hline
 3W4Q &  773 & \textbf{0.737} &          0.725 &          0.649 &          0.628 &  4G6C &  658 & \textbf{0.591} &           0.59 &          0.528 &          0.501 \\ \hline
 3ZBD &  213 &          0.651 &          0.516 &          0.632 & \textbf{0.776} &  4G7X &  194 & \textbf{0.688} &          0.587 &          0.624 &          0.622 \\ \hline
 3ZIT &  152 &  \textbf{0.43} &          0.404 &          0.392 &          0.412 &  4GA2 &  144 &          0.528 &          0.485 &          0.406 & \textbf{0.642} \\ \hline
 3ZRX &  221 &  \textbf{0.59} &          0.562 &          0.391 &          0.501 &  4GMQ &   92 &          0.678 &          0.628 &           0.55 & \textbf{0.862} \\ \hline
 3ZSL &  138 &          0.691 &          0.687 &          0.526 &  \textbf{0.87} &  4GS3 &   90 &          0.544 &          0.522 &          0.547 & \textbf{0.778} \\ \hline
 3ZZP &   74 &          0.524 &           0.46 &          0.448 & \textbf{0.944} &  4H4J &  236 &           0.81 &          0.806 &          0.689 & \textbf{0.832} \\ \hline
 3ZZY &  226 & \textbf{0.746} &          0.709 &          0.728 &          0.621 &  4H89 &  168 & \textbf{0.682} &          0.588 &          0.596 &          0.672 \\ \hline
 4A02 &  166 &          0.618 &          0.516 &          0.303 & \textbf{0.748} &  4HDE &  168 &          0.745 &          0.728 &          0.615 & \textbf{0.911} \\ \hline
 4ACJ &  167 &          0.748 &          0.746 & \textbf{0.759} &          0.752 &  4HJP &  281 &          0.703 &          0.649 &           0.51 & \textbf{0.764} \\ \hline
 4AE7 &  186 & \textbf{0.724} &          0.717 &          0.717 &          0.717 &  4HWM &  117 &          0.638 &          0.622 &          0.499 & \textbf{0.767} \\ \hline
 4AM1 &  345 & \textbf{0.674} &          0.619 &           0.46 &          0.665 &  4IL7 &   85 &          0.446 &          0.404 &          0.316 & \textbf{0.635} \\ \hline
 4ANN &  176 &          0.551 &          0.536 &           0.47 & \textbf{0.641} &  4J11 &  357 &  \textbf{0.62} &          0.562 &          0.401 &          0.528 \\ \hline
 4AVR &  188 &           0.68 &          0.605 &           0.65 & \textbf{0.784} &  4J5O &  220 &          0.793 &          0.757 &          0.777 & \textbf{0.824} \\ \hline
 4AXY &   54 &            0.7 &          0.623 &  \textbf{0.72} &          0.696 &  4J5Q &  146 &          0.742 &          0.742 &          0.689 & \textbf{0.787} \\ \hline
 4B6G &  558 & \textbf{0.765} &          0.756 &          0.669 &          0.758 &  4J78 &  305 & \textbf{0.658} &          0.648 &          0.608 &          0.574 \\ \hline
 4B9G &  292 & \textbf{0.844} &          0.816 &          0.763 &          0.832 &  4JG2 &  185 &          0.746 &          0.736 &          0.543 & \textbf{0.775} \\ \hline
 4DD5 &  387 &          0.615 &          0.596 &          0.351 & \textbf{0.712} &  4JVU &  207 & \textbf{0.723} &          0.697 &          0.553 &          0.697 \\ \hline
 4DKN &  423 & \textbf{0.781} &          0.761 &          0.539 &          0.728 &  4JYP &  534 &          0.688 &          0.682 &          0.538 & \textbf{0.711} \\ \hline
 4DND &   95 & \textbf{0.763} &           0.75 &          0.582 &          0.623 &  4KEF &  133 &           0.58 &           0.53 &          0.324 & \textbf{0.701} \\ \hline
 4DPZ &  109 &           0.73 &          0.726 &          0.651 & \textbf{0.795} &  5CYT &  103 &          0.441 &          0.421 &          0.331 &   \textbf{0.7} \\ \hline
 4DQ7 &  328 &           0.69 &          0.683 &          0.376 & \textbf{0.763} &  6RXN &   45 &          0.614 &          0.574 &          0.594 & \textbf{0.825} \\ \hline
 
 \caption{The comparison of correlation coefficient of mGLI with previous methods including opFRI, prFRI, and GNM. N refers to the number of residues in the protein. And the best value for each protein is marked in bold.}\label{364}
\end{longtable}

	\bibliographystyle{unsrtnat}
	\bibliography{ref}

\end{document}